%%%%%%%%%%%% MACROS FOR PAPERS AND REPORTS %%%%%%%%%%%%
%My little macros

\def\ignore#1{}
 
%%%%%%%%%%%%%%%%%%%%%

\newcount\sectnum
\newcount\subsectnum
\newcount\eqnumber

\global\eqnumber=1\sectnum=0

% Equation labels

\def\lab{(\the\sectnum.\the\eqnumber)}

%Example of use: suppose we want to give a label \lgh to an equation
% $$ ......  \xdef\lgh{\lab} \eqnum \show{lgh}$$
% Later refer to Eq. \lgh\ ...
% Note the \ after \lgh; it seems to be needed if we want the equation number
% to be followed by a space; not needed if followed by . or ,

%The next macro is used to display labels in drafts, so that you do
%not have to remember them

\def\show#1{#1}

%The next macro is to be used for final drafts that do not display labels
%\def\show#1{}

%%%%%%%%%%%%%%%%%%%%%

\def\smskip{\vskip 5 pt}
\def\medskip{\vskip 10 pt}
\def\bigskip{\vskip 15 pt}
\def\pn{\par\noindent}
\def\br{\break}

\def\bl{\bigl} 
\def\br{\bigr}

\def\frac#1#2{{#1\over #2}}

\def\ol#1{\overline{#1}}

\def\e{\epsilon}

 %break line; horizontal space

\def\old#1{}% invalidates text in braces 
\def\leaderfill{\leaders\hbox to 1em{\hss.\hss}\hfill}
%Example of use: \line{1. Optimality Conditions\leaderfill p.\ 2}

% John's macros

\parindent=2pc
\baselineskip=15pt
\vsize=8.7 true in
\voffset=0.125 true in
\parskip=3pt

% vector/matrix macros

%eqalign macros
\def\minprob#1#2#3{$$\eqalign{&\hbox{minimize\ \ }#1\cr &\hbox{subject to\ \
}#2\cr}\ifnum 0=#3{}\else\eqno(#3)\fi$$}        
     
\def\maxprob#1#2#3{$$\eqalign{&\hbox{maximize\ \ }#1\cr &\hbox{subject to\ \
}#2\cr}\ifnum 0=#3{}\else\eqno(#3)\fi$$}        
     
\def\aligntwo#1#2#3#4#5{$$\eqalign{#1&#2\cr #3&#4\cr}
\ifnum 0=#5{}\else\eqno(#5)\fi$$}
\def\alignthree#1#2#3#4#5#6#7{$$\eqalign{#1&#2\cr #3&#4\cr #5&#6\cr}
\ifnum 0=#7{}\else\eqno(#7)\fi$$}

% Macros to automatically advance equation and other numbers

\def\eqnum{\eqno{\hbox{(\the\sectnum.\the\eqnumber)}\global\advance\eqnumber
by1}}

\def\eqnu{\eqno{\hbox{(\the\sectnum.\the\eqnumber)}\global\advance\eqnumber
by1}}

\newcount\examplnumber
\def\examplnum{\global\advance\examplnumber by1}

\newcount\figrnumber
\def\figrnum{\global\advance\figrnumber by1}

\newcount\propnumber
\def\propnum{\global\advance\propnumber by1}

\newcount\defnumber
\def\defnum{\global\advance\defnumber by1}

\newcount\lemmanumber
\def\lemmanum{\global\advance\lemmanumber by1}

\newcount\assumptionnumber
\def\assumptionnum{\global\advance\assumptionnumber by1}

\newcount\conditionnumber
\def\conditionnum{\global\advance\conditionnumber by1}

\def\exampl{\the\sectnum.\the\examplnumber}
\def\figr{\the\sectnum.\the\figrnumber}
\def\propn{\the\sectnum.\the\propnumber}
\def\defn{\the\sectnum.\the\defnumber}
\def\lemman{\the\sectnum.\the\lemmanumber}
\def\assumptionn{\the\sectnum.\the\assumptionnumber}
\def\condn{\the\sectnum.\the\conditionnumber}

\def\section#1{\goodbreak\vskip 3pc plus 6pt minus 3pt\leftskip=-2pc
   \global\advance\sectnum by 1\eqnumber=1
\global\examplnumber=1\figrnumber=1\propnumber=1\defnumber=1\lemmanumber=1\assumptionnumber=1 \conditionnumber =1 \subsectnum=0%
   \line{\hfuzz=1pc{\hbox to 3pc{\bf %\the\sectnum.\quad
   \vtop{\hfuzz=1pc\hsize=38pc\hyphenpenalty=10000\noindent\uppercase{\the\sectnum.\quad #1}}\hss}}
			\hfill}
			\leftskip=0pc\nobreak\tenf
			\vskip 1pc plus 4pt minus 2pt\noindent\ignorespaces}

% ETP Macros

%\def\section#1{\goodbreak\vskip 3pc plus 6pt minus 3pt\leftskip=-2pc
%   \global\advance\sectnum by 1\eqnumber=1
%   \line{\hfuzz=1pc{\hbox to 3pc{\bf %\the\sectnum.\quad
%   \vtop{\hfuzz=1pc\hsize=38pc\hyphenpenalty=10000\noindent\uppercase{#1}}\hss}}
%                        \hfill}
%                        \leftskip=0pc\nobreak\tenf
%                        \vskip 1pc plus 4pt minus 2pt\noindent\ignorespaces}

\def\sect#1{\noindent\leftskip=-2pc\tenf
   \goodbreak\vskip 1pc plus 4pt minus 2pt
                \global\advance\subsectnum by 1\eqnumber=1
   \line{\hfuzz=1pc{\hbox to 3pc{\bf %\the\sectnum.\quad
   \vtop{\hfuzz=1pc\hsize=38pc\hyphenpenalty=10000\noindent\uppercase{{\bf #1}}}\hss}}
                        \hfill}
   \leftskip=0pc\nobreak\tenf
                        \vskip 1pc plus 4pt minus 2pt\nobreak\noindent\ignorespaces}

\def\subsection#1{\noindent%
   \goodbreak\vskip 1pc plus 4pt minus 2pt%
 		\global\advance\subsectnum by 1%
   \line{\hfuzz=1pc{\hbox to 3pc%
   {\bf  \vtop{\hfuzz=1pc\hsize=38pc\hyphenpenalty=10000\noindent{\bf 
  \the\sectnum.\the\subsectnum\ \ \ #1}}\hss}}%
			\hfill}%
   \nobreak%
			\vskip 1pc plus 4pt minus 2pt\nobreak\noindent\ignorespaces}%

\def\subsubsection#1{\goodbreak\vskip 1pc plus 4pt minus 2pt
   \hfuzz=3pc\leftskip=0pc\noindent\tenit #1 \nobreak\tenf\vskip 6pt plus 1pt
                                minus 1pt\nobreak\ignorespaces\leftskip=0pc}
%
%\def\rthl{Sec. \the\chapnum.\the\sectnum}                      
%\def\rthc{#1}\nobreak\noindent\ignorespaces
%\newcount\sectnum \sectnum=0
%\newcount\subsectnum \subsectnum=0

\def\beginexample#1{\noindent\goodbreak\vskip 6pt plus 1pt minus 1pt
\noindent
  \hbox {\bf Example #1\hss}%\break%\noindent
  \nobreak\vskip 4pt plus 1pt minus 1pt \nobreak\noindent\ninef
  \global\advance
                \leftskip by\parindent\pn}
\def\endexample{\vskip 12pt\tenf\par
  \global\advance\leftskip by -\parindent
  }

\def\beginexercise#1{\noindent\goodbreak\vskip 6pt plus 1pt minus 1pt \noindent\global\normalbaselineskip=12pt
  \hbox {\bf Exercise #1\hss}%\break%\noindent
  \nobreak\vskip 4pt plus 1pt minus 1pt 
  \nobreak\noindent\ninef\global\advance\leftskip
                        by\parindent\pn}
\def\endexercise{\vskip 12pt\tenf\par
  \global\advance\leftskip by -\parindent
  }

\def\beginsection#1{\noindent\goodbreak\vskip 6pt plus 1pt minus 1pt \noindent\global\normalbaselineskip=12pt
  \hbox {\it #1\hss}
  \vskip 0.1pt plus 1pt minus 1pt \nobreak\noindent\ninef\global\advance
                \leftskip by\parindent\noindent\pn}
\def\endsection{\vskip 12pt\tenf\par
  \global\advance\leftskip by -\parindent
}

%

% Header/Title macros

\def\proposition#1{\smskip\pn{\bf Proposition #1}\quad}

\def\ref{\smskip\pn}

\def\chapter#1#2{{\bf \centerline{\helbigbig
{#1}}}\bigskip\bigskip{\bf \centerline{\helbigbig
{#2}}}\bigskip\bigskip} % ex. \chapter{Chapter 1}{Title of chapter}

 % ex. \longchapter{Chapter 1}{Title of chapter}{Title of
 %chapter}

 % ex. \papertitle{Title of paper}{Names of Authors}

\def\longpapertitle#1#2#3{{\bf \centerline{\helbigb
{#1}}}\bigskip{\bf \centerline{\helbigb
{#2}}}\bigskip\bigskip{\centerline{
by}}\bigskip{\bf \centerline{
{#3}}}\bigskip\bigskip} 
% ex. \longpapertitle{First part of title of paper}
%{2nd part of title of paper}{Names of Authors}

% List macros

\def\nitem#1{\smskip\item{#1}}

\newcount\alphanum
\newcount\romnum

\def\alphaenumerate{\ifcase\alphanum \or (a)\or (b)\or (c)\or (d)\or (e)\or
(f)\or (g)\or (h)\or (i)\or (j)\or (k)\fi}
\def\romenumerate{\ifcase\romnum \or (i)\or (ii)\or (iii)\or (iv)\or (v)\or
(vi)\or (vii)\or (viii)\or (ix)\or (x)\or (xi)\fi}

\def\alist{\begingroup\vskip10pt\alphanum=1% alphabetical list
\parskip=2pt\parindent=0pt \leftskip=3pc
\everypar{\llap{{\rm\alphaenumerate\hskip1em}}\advance\alphanum by1}}

\def\nolist{\begingroup\vskip10pt\alphanum=0% numerical list
\parskip=2pt\parindent=0pt \leftskip=3pc
\everypar{\llap{\global\advance\alphanum by1(\the\alphanum)\hskip1em}}}

\def\romlist{\begingroup\vskip10pt\romnum=1% roman list
\parskip=2pt\parindent=0pt \leftskip=5pc
\everypar{\llap{{\rm\romenumerate\hskip1em}}\advance\romnum by1}}

% romlist indents more than alist or nolist and can be used inside them

%Figure, table, and box macros

\long\def\fig#1#2#3{\vbox{\vskip1pc\vskip#1
\prevdepth=12pt \baselineskip=12pt
\vskip1pc
\hbox to\hsize{\hfill\vtop{\hsize=25pc\noindent{\eightbf Figure #2\ }
{\eightpoint#3}}\hfill}}}%Figure space definition. Example of use:
%\fig{16pc}{1.1}{A network with one central processor and a separate
%communication link to each device.}

\long\def\widefig#1#2#3{\vbox{\vskip1pc\vskip#1
\prevdepth=12pt \baselineskip=12pt
\vskip1pc
\hbox to\hsize{\hfill\vtop{\hsize=28pc\noindent{\eightbf Figure #2\ }
{\eightpoint#3}}\hfill}}}

\long\def\table#1#2{\vbox{\vskip0.5pc
\prevdepth=12pt \baselineskip=12pt
\hbox to\hsize{\hfill\vtop{\hsize=25pc\noindent{\eightbf Table #1\ }
{\eightpoint#2}}\hfill}}}

%Running Head Macros
 
\def\rightheadline#1{\headline{\tenrm\hfil #1}}

% Concept Macros

\long\def\leftfig#1#2{\vbox{\smskip\hsize=220pt
\vtop{{\noindent {\bf #1}}}
\smskip
\noindent
\vbox{{\noindent #2}}
}}

\long\def\rightfig#1#2#3{\vbox{\smskip\vskip#1
\prevdepth=12pt \baselineskip=12pt
\hsize=210pt
\smskip
\vbox{\noindent{\eightbold #2}
\hskip1em{\eightpoint#3}}
}}

\long\def\concept#1#2#3#4#5{\bigskip\hrule
\vbox{\hbox{\leftfig{#1}{#2} \hskip3em
\rightfig{#3}{#4}{#5}} \smskip}
\hrule\bigskip}

% Example of Use: \concept{Title of Concept}{Text}
% {Figure size}{Figure number?}{Figure caption}

\long\def\bconcept#1#2#3#4#5#6#7{
\vbox{
\hbox to \hsize{\vtop{\par #1}}
\concept{#2}{#3}{#4}{#5}{#6}
\hbox to \hsize{\vtop{\par #7}}
\smskip}
}

% Example of Use: \bconcept{Preceding text}{Title of Concept}{Text}
% {Figure size}{Figure number}{Figure caption}{Following text}

% same as concept but without the \hrule's; ready to be boxed

% Put inside a box

\def\boxit#1{\vbox{\hrule\hbox{\vrule\kern3pt
                                \vbox{\kern3pt#1\kern3pt}\kern3pt\vrule}\hrule}}
% example of use: \setbox0=\vbox{.... }; \boxit{\box0}
\def\centerboxit#1{$$\vbox{\hrule\hbox{\vrule\kern3pt
                                \vbox{\kern3pt#1\kern3pt}\kern3pt\vrule}\hrule}$$}
% example of use: \setbox0=\vbox{.... }; \centerboxit{\box0}

\long\def\boxtext#1#2{$$\boxit{\vbox{\hsize #1\noindent\strut #2\strut}}$$}
% example of use: \boxtext{462pt}{This is the boxed text.}; 462pt is max length

% Picture macros and examples
%
% figures must be pasted from mcdraw
%
% Look in the 'Windows' menu for the pictures window
% It's like the Scrapbook -- cut and paste pictures
%

\def\picture #1 by #2 (#3){
  \vbox to #2{
    \hrule width #1 height 0pt depth 0pt
    \vfill
    \special{picture #3} % this is the low-level interface
    }
  }
% The first dimension of the picture macro is the width the second is depth

\def\scaledpicture #1 by #2 (#3 scaled #4){{
  \dimen0=#1 \dimen1=#2
  \divide\dimen0 by 1000 \multiply\dimen0 by #4
  \divide\dimen1 by 1000 \multiply\dimen1 by #4
  \picture \dimen0 by \dimen1 (#3 scaled #4)}
  }

%
% Note that you can also say, e.g.,
%  \special{postscript xxx yyy zzz}
% to include PostScript graphics in your documents
%
%Examples of use
%\def\stripes{\picture 2.29in by 1.75in (AWstripes)}
% By executing \stripes 
%\def\annie{\scaledpicture 102pt by 239pt (annie scaled 2000)}
%\def\finder{\picture 260pt by 165pt (screen0 scaled 500)}
%\def\icon{\picture 7in by 7in (icon)}
%Example of use
%\annie
%Example of centered picture \line{\hfil\annie\hfil}

%Figure w/  caption macro
\long\def\captfig#1#2#3#4#5{\vbox{\vskip1pc
\hbox to\hsize{\hfill{\picture #1 by #2 (#3)}\hfill}
\prevdepth=9pt \baselineskip=9pt
\vskip1pc
\hbox to\hsize{\hfill\vtop{\hsize=24pc\noindent{\eightbold Figure #4}
\hskip1em{\eightpoint#5}}\hfill}}}

%Examples of use of Figure macros
%\captfig{8.53pc}{19.9pc}{picturename}{5}{caption.}
%\captfig{2.23in}{2in}{picturename scaled 500}{16}{Caption.}
%The macro centers the picture.
%The first two numbers should be the true width
% and height after the picture has been scaled.
% So if the picture is scaled by 50% (500), the width and height in
% the macro should onw half of what they would be if the picture
% is not scaled (1000).
%
%
%
% Postcript macros

\def\illustration #1 by #2 (#3){
  \vskip#2\hskip#1\special{illustration #3} % this is the low-level interface
    }

\def\scaledillustration #1 by #2 (#3 scaled #4){{
  \dimen0=#1 \dimen1=#2
  \divide\dimen0 by 1000 \multiply\dimen0 by #4
  \divide\dimen1 by 1000 \multiply\dimen1 by #4
  \illustration \dimen0 by \dimen1 (#3 scaled #4)}
  }

% SHADEBOX.BSR MACROS
% Author: Leo@vaxc.cc.monash.edu.au
% Original Source:  Posted by Jimm Herreron <HERRERON@SMCVAX.BITNET>
% Modified from the file SHADEBOX.TEX on 9/30/93 by Becky Kaluza of Blue Sky
% Research to work with Textures 1.5 or later.

\newbox\graybox
\newdimen\xgrayspace
\newdimen\ygrayspace
%
% This macro can be used to typeset some text in a framed box with a
% shaded background. A set of examples can be found at the end of this
% file.
%
% This is a plain \TeX\ file modified for use on the Macintosh with Textures
% 1.5 or later.
%
% The characteristics of the shaded boxes are controlled by the following
% parameters
%
%   \xgrayspace = the space added before and after the text
%   \ygrayspace = the space above and below the text
%   \grayshade  = the gray colour 0 = black 1 = white
%   \linewidth  = the thickness of the border in points
%   \theradius  = the radius of the rounded corners in points
%   \thevskip   = extra \vskip added above and below the shaded box
%                 (applies only to \parashade)
%
%----------------------------------------------------------------------------
%
% The following \TeX code was adapted from previous work by
%
%            Je'ro^me Maillot, maillot@bora.inria.fr
%----------------------------------------------------------------------------
%
% Use the following for one or more words within a line.
%

\def\Textshade#1#2#3#4#5#6{%
    \xgrayspace=#4pt%
    \ygrayspace=#4pt%
    \def\grayshade{#3}%
    \def\linewidth{#5}%
    \def\theradius{#6}%
    \setbox\graybox=\hbox{\surroundboxa{#2}}%
    \hbox{%
    \hbox to 0pt{%
%!    \special{"gsave newpath 0 0 moveto                                %
    \PScommands
}%
    \box\graybox}}%
%
% Use the following for paragraphs.
%
\long%

\long%
\def\Parashade#1#2#3#4#5#6#7{%
    \xgrayspace=#4pt%
    \ygrayspace=#4pt%
    \def\grayshade{#3}%
    \def\linewidth{#5}%
    \def\theradius{#6}%
    \def\thevskip{#7pt}%
    \setbox\graybox=\hbox{\surroundboxb{#2}}%
    \vskip\thevskip%
    \hbox{%
    \hbox to 0pt{%
%!    \special{"gsave newpath 0 0 moveto                                %
    \PScommands
}%
     \box\graybox}%
     \vskip\thevskip%
}%
%----------------------------------------------------------------------------
%
% A pair of box macros. Each builds a slightly oversized box
% containing the text. The text is centred both in the vertical
% horizontal directions.
%
% Use the following for one or more words within a line.
%
\long\def\surroundboxa#1{\leavevmode\hbox{\vtop{%
\vbox{\kern\ygrayspace%
\hbox{\kern\xgrayspace#1%
      \kern\xgrayspace}}\kern\ygrayspace}}}
%
% Use the following for a paragraphs.
%
\long\def\surroundboxb#1{\leavevmode\hbox{\vtop{%
\vbox{\kern\ygrayspace%
\hbox{\kern\xgrayspace\vbox{\advance\hsize-2\xgrayspace#1}%
      \kern\xgrayspace}}\kern\ygrayspace}}}
%----------------------------------------------------------------------------
%
% Here are some simple PostScript routines.
%
% The TeX command \PScommands must be called before any of the
% shading routines can be used.
%
%!\long\def\PScommands{\special{! TeXDict begin
\long\def\PScommands{%
\special{rawpostscript
/sharpbox{%
           newpath
           xmin ymin moveto
           xmin ymax lineto
           xmax ymax lineto
           xmax ymin lineto
           xmin ymin lineto
           closepath 
          }bind def
}%
\special{rawpostscript
/sharpboxnb{%
           newpath
           xmin ymin moveto
           xmin ymax lineto
           xmax ymax lineto
           xmax ymin lineto
%           xmin ymin lineto
%           closepath 
          }bind def
}%
\special{rawpostscript
/sharpboxnt{%
           newpath
           xmin ymax moveto
           xmin ymin lineto
           xmax ymin lineto
           xmax ymax lineto
%           xmin ymin lineto
%           closepath 
          }bind def
}%
\special{rawpostscript
/roundbox{%
           newpath
           xmin radius add ymin moveto
           xmax ymin xmax ymax radius arcto
           xmax ymax xmin ymax radius arcto
           xmin ymax xmin ymin radius arcto
           xmin ymin xmax ymin radius arcto 16 {pop} repeat
           closepath
          }bind def
}%
\special{rawpostscript
/sharpcorners{%
               sharpbox gsave grayshade setgray fill grestore 
               linewidth setlinewidth stroke
              }bind def
}%
\special{rawpostscript
/sharpcornersnt{%
               sharpboxnt gsave grayshade setgray fill grestore 
               linewidth setlinewidth stroke
              }bind def
}%
\special{rawpostscript
/sharpcornersnb{%
               sharpboxnb gsave grayshade setgray fill grestore 
               linewidth setlinewidth stroke
              }bind def
}%
\special{rawpostscript
/roundcorners{%
               roundbox gsave grayshade setgray fill grestore 
               linewidth setlinewidth stroke
              }bind def
}%
\special{rawpostscript
/plainbox{%
           sharpbox grayshade setgray fill 
          }bind def
}%
% Here are the two new options
%
\special{rawpostscript
/roundnoframe{%
               roundbox grayshade setgray fill 
              }bind def
}%
\special{rawpostscript
/sharpnoframe{%
               sharpbox grayshade setgray fill 
              }bind def
}%
%!end}%
}%
%
% The \PScommands macro must be invoked before the shaded box macros.
%
%!\PScommands
% To use this, type \textshade{plainbox} or \textshade{roundbox} or
% \textshade{sharpbox}

\def\pshade#1{\Parashade{sharpcorners}{#1}{0.95}{10}{0.5}{10}{10}}

%%%%% BOXES FOR TEXSHOP %%%%%

\def\boxit#1{\vbox{\hrule\hbox{\vrule\kern3pt
                                \vbox{\kern3pt#1\kern3pt}\kern3pt\vrule}\hrule}}
% example of use: \setbox0=\vbox{.... }; \boxit{\box0}

\def\boxitnb#1{\vbox{\hrule\hbox{\vrule\kern3pt
                                \vbox{\kern3pt#1\kern3pt}\kern3pt\vrule}}}

\def\boxitnt#1{\vbox{\hbox{\vrule\kern3pt
                                \vbox{\kern3pt#1\kern3pt}\kern3pt\vrule}\hrule}}

\def\boxitntnb#1{\vbox{\hbox{\vrule\kern3pt
                                \vbox{\kern3pt#1\kern3pt}\kern3pt\vrule}}}

\long\def\boxtext#1#2{$$\boxit{\vbox{\hsize #1\noindent\strut #2\strut}}$$}
% example of use: \boxtext{462pt}{This is the boxed text.}; 462pt is max length
\long\def\boxtextnb#1#2{$$\boxitnb{\vbox{\hsize #1\noindent\strut #2\strut}}$$}
% example of use: \boxtext{462pt}{This is the boxed text.}; 462pt is max length
\long\def\boxtextnt#1#2{$$\boxitnt{\vbox{\hsize #1\noindent\strut #2\strut}}$$}
% example of use: \boxtext{462pt}{This is the boxed text.}; 462pt is max length

% example of use: \boxtext{462pt}{This is the boxed text.}; 462pt is max length

\def\texshopbox#1{\boxtext{462pt}{\vskip-1.5pc\pshade{\vskip-1.0pc#1\vskip-2.0pc}}}
\def\texshopboxnt#1{\boxtextnt{462pt}{\vskip-1.5pc\pshade{\vskip-1.0pc#1\vskip-2.0pc}}}
\def\texshopboxnb#1{\boxtextnb{462pt}{\vskip-1.5pc\pshade{\vskip-1.0pc#1\vskip-2.0pc}}}

%***************************************************
%         FONTS
%***************************************************

% ROMAN
%
%
%
%
%
%
%
%
\font\helbigbig=cmr10 scaled 2500%
\font\helbigb=cmbx10 scaled 1500%
\font\eightbold=cmbx8%

\def\tenf{\hel}%
\def\tenit{\heli}%
\def\ninef{\ninehel}%
\def\nineit{\nineheli}%
%
%

%  FONT FAMILIES

\font\tenrm=cmr10%
\font\teni=cmmi10%
\font\tensy=cmsy10%
\font\tenbf=cmbx10%
\font\tentt=cmtt10%
\font\tenit=cmti10%
\font\tensl=cmsl10%

\def\tenpoint{\def\rm{\fam0\tenrm}%
\textfont0=\tenrm%
\textfont1=\teni%
\textfont2=\tensy%
\textfont\itfam=\tenit%
\textfont\slfam=\tensl%
\textfont\ttfam=\tentt%
\textfont\bffam=\tenbf%
\scriptfont0=\sevenrm%
\scriptfont1=\seveni%
\scriptfont2=\sevensy%
%\scriptfont3=\tenex%
\scriptscriptfont0=\sixrm%
\scriptscriptfont1=\sixi%
\scriptscriptfont2=\sixsy%
%\scriptscriptfont3=\tenex%
\def\it{\fam\itfam\tenit}%
\def\tt{\fam\ttfam\tentt}%
\def\sl{\fam\slfam\tensl}%
\scriptfont\bffam=\sevenbf%
\scriptscriptfont\bffam=\sixbf%
\def\bf{\fam\bffam\tenbf}%
\normalbaselineskip=18pt%
\normalbaselines\rm}%

\font\ninerm=cmr9%
\font\ninebf=cmbx9%
\font\nineit=cmti9%
\font\ninesy=cmsy9%
\font\ninei=cmmi9%
\font\ninett=cmtt9%
\font\ninesl=cmsl9%

\def\ninepoint{\def\rm{\fam0\ninerm}%
\textfont0=\ninerm%
\textfont1=\ninei%
\textfont2=\ninesy%
\textfont\itfam=\nineit%
\textfont\slfam=\ninesl%
\textfont\ttfam=\ninett%
\textfont\bffam=\ninebf%
\scriptfont0=\sixrm%
\scriptfont1=\sixi%
\scriptfont2=\sixsy%
%\scriptfont3=\tenex%
\def\it{\fam\itfam\nineit}%
\def\tt{\fam\ttfam\ninett}%
\def\sl{\fam\slfam\ninesl}%
\scriptfont\bffam=\sixbf%
\scriptscriptfont\bffam=\fivebf%
\def\bf{\fam\bffam\ninebf}%
\normalbaselineskip=16pt%
\normalbaselines\rm}%

\font\eightrm=cmr8%
\font\eighti=cmmi8%
\font\eightsy=cmsy8%
\font\eightbf=cmbx8%
\font\eighttt=cmtt8%
\font\eightit=cmti8%
\font\eightsl=cmsl8%

\def\eightpoint{\def\rm{\fam0\eightrm}%
\textfont0=\eightrm%
\textfont1=\eighti%
\textfont2=\eightsy%
\textfont\itfam=\eightit%
\textfont\slfam=\eightsl%
\textfont\ttfam=\eighttt%
\textfont\bffam=\eightbf%
\scriptfont0=\sixrm%
\scriptfont1=\sixi%
\scriptfont2=\sixsy%
%\scriptfont3=\tenex%
\scriptscriptfont0=\fiverm%
\scriptscriptfont1=\fivei%
\scriptscriptfont2=\fivesy%
%\scriptscriptfont3=\tenex%
\def\it{\fam\itfam\eightit}%
\def\tt{\fam\ttfam\eighttt}%
\def\sl{\fam\slfam\eightsl}%
%\scriptfont\bffam=\sixbf%
\scriptscriptfont\bffam=\fivebf%
\def\bf{\fam\bffam\eightbf}%
\normalbaselineskip=14pt%
\normalbaselines\rm}%

\font\sevenrm=cmr7%
\font\seveni=cmmi7%
\font\sevensy=cmsy7%
\font\sevenbf=cmbx7%

\def\sevenpoint{%
   \def\rm{\sevenrm}\def\bf{\sevenbf}%
   \def\smc{\sevensmc}\baselineskip=12pt\rm}%

\font\sixrm=cmr6%
\font\sixi=cmmi6%
\font\sixsy=cmsy6%
\font\sixbf=cmbx6%

\fontdimen13\tensy=2.6pt%
\fontdimen14\tensy=2.6pt%
\fontdimen15\tensy=2.6pt%
\fontdimen16\tensy=1.2pt%
\fontdimen17\tensy=1.2pt%
\fontdimen18\tensy=1.2pt%       

\def\tenf{\tenpoint}%
\def\ninef{\ninepoint}%
%

%%%%%%%%%%%% END OF MACROS %%%%%%%%%%%%

%%%%%%%%%% REDEFINITION OF BOX SPACING %%%%%%%%%%%%%%%%

\def\texshopbox#1{\boxtext{462pt}{\vskip-1.5pc\pshade{\vskip-1.0pc#1\vskip-2.0pc}}}
\def\texshopboxnt#1{\boxtextnt{462pt}{\vskip-1.5pc\pshade{\vskip-1.0pc#1\vskip-2.0pc}}}
\def\texshopboxnb#1{\boxtextnb{462pt}{\vskip-1.5pc\pshade{\vskip-1.0pc#1\vskip-2.0pc}}}

\long\def\fig#1#2#3{\vbox{\vskip1pc\vskip#1
\prevdepth=12pt \baselineskip=12pt
\vskip1pc
\hbox to\hsize{\hfill\vtop{\hsize=30pc\noindent{\eightbf Figure #2\ }
{\eightpoint#3}}\hfill}}}

\def\show#1{}
\def\frac#1#2{{#1\over #2}}
\def\jstar{J^{\raise0.04pt\hbox{\sevenpoint *}} }

\rightheadline{\botmark}

\pageno=1

%%%%%%%%%%%%%%%%%%%%%%%%%%%%%%%%%%%%%%%%%%%%

% GRAPHICS MACROS FOR TEX

%\input graphicx.texomplentary slackness

% DO NOT USE THE LATEX CODE

%\input miniltx
%
%\ifx\pdfoutput\undefined
 % \def\Gin@driver{dvips.def} % we are not running PDFTeX
%\else
%  \def\Gin@driver{pdftex.def} % we are running PDFTeX
%\fi

%\input graphicx.sty
%\resetatcatcode
%

% ALTERNATIVE SIMPLE METHOD

\input epsf

% Example of Use

%\epsfbox{PI_Minimax.eps}
%\centerline{\hskip0pc\epsfxsize = 4.5in \epsfbox{PI_Minimax.eps}}
%\centerline{\epsfscale=800 \epsfbox{PI_Minimax.eps}} NOTE: \epsfscale is questionable

%%%%%%%%%%%%%%%%%%%%%%%%%%%%%%%%%%%%%%%%%%%%

\pn {\bf July 2022}\hfill{\bf Arizona State University/SCAI Report}
%\hfill{\bf Report LIDS 2646.}
%\hfill {\bf Submitted for publication to JOTA}%
\bigskip \bigskip

\bigskip\bigskip\bigskip

\def\longpapertitle#1#2#3{{\bf \centerline{\helbigb
{#1}}}\medskip{\bf \centerline{\helbigb
{#2}}}\medskip{\centerline{
by}}\medskip{\bf \centerline{
{#3}}}\bigskip}

\longpapertitle{New Auction Algorithms for Path Planning,}{Network Transport, and Reinforcement Learning}{{Dimitri Bertsekas\footnote{\dag}{\ninepoint Fulton Professor of Computational Decision Making, School of Computing and Augmented Intelligence, Arizona State University, Tempe, AZ.}}}

\centerline{\bf Abstract}

We consider some classical optimization problems in path planning and network transport, and we introduce new auction-based algorithms for their optimal and suboptimal solution. The algorithms are based on mathematical ideas that are related to competitive bidding by persons for objects and the attendant market equilibrium, which underlie auction processes. However, the starting point  of our algorithms is different, namely weighted and unweighted path construction in directed graphs, rather than assignment of persons to objects. The new algorithms have several potential advantages over existing methods: they are empirically faster in some important contexts, such as max-flow, they are well-suited for on-line replanning, and they can be adapted to distributed asynchronous operation. Moreover, they allow arbitrary initial prices, without complementary slackness restrictions, and thus are better-suited to take advantage of reinforcement learning methods that use off-line training with data, as well as on-line training during real-time operation. The new algorithms may also find use in reinforcement learning contexts involving approximation, such as multistep lookahead and tree search schemes, and/or rollout algorithms.

\vfill\eject

\section{Introduction}
%\vskip-1pc

\pn In this paper, we introduce new auction algorithms for broad classes of network flow problems. Our proposed methodology aims to improve the efficiency and flexibility of existing auction algorithms for linear single commodity network optimization,  including shortest path planning, matching, assignment, and network transportation problems. 

Auction algorithms for network flow optimization have a long history, starting with the author's original paper [Ber79] that dealt with the assignment problem. They are discussed in detail in many sources, including the author's book [Ber98] and tutorial survey [Ber92] (which contain many references, too numerous to list here).  Besides excellent computational complexity properties, their advantages over other network flow methods include their suitability for reoptimization and parallelization.  They have found use in contexts involving classical weighted matching, route planning, scheduling, and related network optimization problems. 

Auction algorithms have also been considered widely in applications of optimal transport, which are currently very popular; see e.g., Brenier et al.\ [BFH03], Villani [Vil09], [Vil21],  Santambrogio [San15], Galichon [Gal16], Schmitzer [Sch16], [Sch19],  Walsh and Dieci [WaD17], [WaD19], Peyre and Cuturi [PeC19], Merigot and Thibert [MeT21], and the references quoted there. Several code implementations of auction algorithms have become publicly available, some of which are accessible from the author's website. 

In this paper, we are also aiming at applications in machine learning, data mining, and artificial intelligence, taking advantage of the ability of our algorithms to adapt to changing environments, and their suitability for off-line and on-line training with data, through the use of machine learning and reinforcement learning techniques. For examples of related contexts, which have involved the use of auction algorithms, see the papers by Kosowsky and Yuille [KoY94], Bayati,  Shah, and Sharma [BSS08], Wang  and Xia [WaX12], Lewis et al.\ [LBD21], Bicciato and Torsello ]BiT22], and Clark et al.\ [CLG22].

The starting point for our development is an algorithm for computing some (not necessarily shortest) path connecting an origin node to a destination node in a directed graph. The algorithm uses node prices to guide the search for the path. A special case of this algorithm was given in the author's paper [Ber95a] as part of an auction algorithm for  max-flow, and was also described in the book [Ber98], Section 3.3.1. The more general version proposed here allows unrestricted choice of the initial prices, which among others, facilitates its use in reinforcement learning and on-line replanning contexts.

Our algorithm also resembles an auction/shortest path algorithm that was proposed in the author's paper [Ber91]. Contrary to the earlier algorithm, the new algorithm uses a positive $\e$ parameter to effect larger price changes, which in turn speed up its convergence. It produces a shortest path for sufficiently small $\e$, and it is also well-suited for an $\e$-scaling approach, a central technique in auction algorithms, which improves computational efficiency.

Our  new path construction algorithm can  be used as a building block for methods to address a host of optimal network flow problems, including matching, assignment, transportation, and minimum cost flow problems with linear as well as separable convex cost functions, in the spirit of the primal-dual/sequential shortest path approach. Auction algorithms for these problems are already available, but our  algorithm is often more flexible. In this paper, we describe briefly the extensions of our algorithm to more general contexts, and we will provide more detailed discussions in future reports and the forthcoming monograph [Ber22a].

 In the next section we will consider a feasibility/path construction  problem, involving a directed graph with a special node $s$, referred to as the {\it origin\/}, and a node $t$, referred to as the {\it destination\/}. We assume that there is at least one directed path from $s$ to $t$, and we want to find one such path.  
One possibility is to use an algorithm that finds a shortest path from $s$ to $t$, with respect to some set of arc lengths or weights; for example a length equal to 1 for every arc. There are several well-known algorithms to solve such a shortest path problem. However, in Section 2 we will introduce a simpler and  faster algorithm that finds a path  from $s$ to $t$ without regard to its optimality. 

In Section 3, we will discuss an extension of the path construction algorithm of Section 2, which uses arc lengths and produces a path that is ``nearly shortest" with respect to these lengths. The algorithms of Sections 2 and 3 bear similarity to an auction algorithm from 1991, which is guaranteed to find a shortest path. This algorithm was proposed in the author's paper [Ber91], and was further described in the book [Ber98], Section 2.6. An important difference is that the earlier algorithm has nonpolynomial complexity,\footnote{\dag}{\ninepoint A polynomial variant of the 
1991 auction algorithm, given in the paper [BPS95], performs very well for single origin-single destination problems with nonnegative arc lengths, but includes features that detract from the flexibility of our new algorithm of Section 3.}
whereas the algorithm of the present paper is polynomial thanks to the use of the $\e$-scaling technique that we noted earlier. Finally, in Sections 4 and 5, we will discuss briefly extensions of our path construction algorithms and their uses in solving matching, assignment, max-flow, transportation, and minimum cost flow problems.
\vskip-1pc

\section{Path Construction in a Directed Graph}

\pn We will now introduce our auction algorithm for path construction.  The algorithm finds just {\it some} path from origin to destination, without aiming for any kind of optimality properties.
It maintains a
path starting at the origin, which at each iteration, is either {\it  extended} by adding a
new node, or {\it contracted} by deleting its terminal node. The decision to extend or contract is based on a  set of variables, one for each node, which are called {\it prices\/}. Roughly speaking, the price of a node is viewed as a measure of the desirability of revisiting and advancing from that node in the future (low-price nodes are viewed as more desirable). Once the destination becomes the terminal node of the path, the algorithm terminates.

\old{
A special case of our algorithm, with more restrictive path and price updating rules, was proposed in the paper [Ber95a] and the book [Ber98], Section 3.3, in connection with the construction of augmenting paths in the context of max-flow problems. The new algorithm to be presented here is more flexible, allowing arbitrary initial prices, and is better suited for on-line replanning.
}

To get an intuitive sense of the algorithm, think of a mouse moving in a graph-like maze, trying
to reach the destination. The mouse criss-crosses the maze, either advancing or backtracking
along its current path, guided by prices that encode how desirable the maze nodes/crosspoints are, based on the mouse's ``learned" experience. The mouse advances forward from high price to low price nodes, going from a node to a downstream neighbor node only if that neighbor has lower price (or equal price under some conditions). It backtracks when it reaches a node whose downstream neighbors all have higher price. In this case, it also suitably increases the price of that node, thus marking the node as less desirable for future exploration, and providing an incentive to explore alternative paths to the destination. 

Our algorithm emulates efficiently the search process just described, guided by a suitable set of rules for price updating. An important side benefit is that the prices provide the means to ``transfer knowledge," in the sense that good learned prices from previous searches can be used as initial prices for subsequent related searches, with an attendant computational speedup. 

\vskip-0.5pc

\subsection{Auction Path Construction Algorithm}
\vskip-0.5pc

\pn We will now describe formally our algorithm, which we call {\it auction/path construction} (APC for short) for finding a path from the origin node $s$ to the destination node $t$ in a directed graph. The arcs of the graph are denoted by $(i,j)$, where $i$ and $j$ are referred to as the {\it start} and {\it end} nodes of the arc. The sets of nodes and the set of arcs are denoted by ${\cal N}$ and ${\cal A}$, respectively. If $(i,j)$ is an arc, it is possible that $(j,i)$ is also an arc. No self arcs of the form $(i,i)$ are allowed. We assume that for any two nodes $i$ and $j$, there is at most one arc with start node $i$ and end node $j$. For any node $i$ we say that node $j$ is a {\it downstream neighbor of $i$} if $(i,j)$ is an arc. 
A node $i$ is called {\it deadend} if it has no downstream neighbors. Note that $s$ is not deadend, since we have assumed that there is a path from $s$ to $t$.

Our algorithm maintains and updates a scalar price $p_i$ for each node $i$. We say that under the current set of prices an arc $(i,j)$ is:
\nitem{(a)} {\it Downhill\/}: If $p_i>p_j$.
\nitem{(b)} {\it Level\/}: If $p_i=p_j$.
\nitem{(c)} {\it Uphill\/}: If $p_i<p_j$.
\smskip
\pn Our algorithm also maintains and updates a directed path $P=(s,n_1,\ldots,n_k)$ that  starts at the origin, and is such that $(s,n_1)$,  and $(n_\ell,n_{\ell+1})$ with $\ell=1,\ldots,k-1$, are arcs. The path is either the degenerate path $P=(s)$, or it ends at some node $n_k\ne s$, which is called the {\it terminal node} of $P$. If $P=(s)$, we also say that the terminal node of $P$ is $s$.

\old{\vskip-0.5pc
\subsection{Algorithm Description}
\vskip-0.5pc
|}

Each iteration of the algorithm starts with a path and a price for each node, which are updated during the iteration using rules that we will now describe. The algorithm starts with the degenerate path $P=(s)$, and with some initial prices, which are arbitrary.\footnote{$\,$\dag}{\ninepoint The arbitrary nature of these prices is a major difference of our algorithm from the earlier auction/path construction algorithms given in [Ber98], Section 2.6 and  3.3. Allowing arbitrary initial prices allows more flexibility in reusing prices from solution of one path finding problem to another similar problem. It also facilitates the use of ``learned" prices that are favorable in similar problem contexts. This property can be important for computational efficiency in many applications.}  It terminates when a path has been found from $s$ to $t$.

At each iteration when the algorithm starts  with a path of the nondegenerate form $P=(s,n_1,\ldots,n_k)$, it either removes from $P$ the terminal node $n_k$ to obtain the new path $\overline P=(s,n_1,\ldots,n_{k-1})$, or it adds to $P$ a node $n_{k+1}$  to obtain the new path $\ol P=(s,n_1,\ldots,n_k,n_{k+1})$. In the former case the operation is called a {\it contraction} to $n_{k-1}$, and in the latter case it is called an {\it extension} to $n_{k+1}$.

At any one iteration the algorithm starts with a path $P$ and a price $p_i$ for each node $i$. At the end of the iteration a new path $\overline P$ is obtained from $P$ through a contraction or an extension. Also the price of the terminal node of $P$ [or the price $p_s$ if $P=(s)$] is increased by a certain amount when there is a contraction. For iterations where the algorithm starts with the degenerate path $P=(s)$, only an extension is possible, i.e., $P=(s)$ is replaced by a path of the form $\overline P=(s,n_1)$.

A key feature of the algorithm, which in fact motivates its design, is that $P$ and the prices $p_i$ satisfy the following property at the start of each iteration for which $P\ne (s)$.

\texshopbox{\pn {\bf Downhill Path Property:}
\pn For the path $P=(s,n_1,\ldots,n_k)$ mantained by the algorithm, the last arc $(n_{k-1},n_k)$ of $P$ is level following a contraction to $n_k$ and is downhill following an extension to $n_k$. All other arcs of $P$ are level.}

The significance of the downhill path property is that {\it when an extension occurs, a cycle cannot be created\/}, in the sense  that the terminal node $n_k$ is different than all the predecessor nodes $s,n_1,\ldots,n_{k-1}$ on the path $P$. The reason is that the downhill path property implies that following an extension, we must have
$$p_{n_k}<p_{n_{k-1}}=p_{n_{k-2}}=\cdots=p_{n_1}=p_s,$$
showing that the terminal node $n_k$ following an extension cannot be equal to any of the preceding nodes of $P$.

In addition to maintaining  the downhill path property, the algorithm is structured so that following a contraction, which changes a nondegenerate path of the form $P=(s,n_1,\ldots,n_k)$ to $\overline P=(s,n_1,\ldots,n_{k-1})$, the price of $n_k$ is increased by a positive amount. In conjunction with the fact that $P$ never contains a cycle, this implies that either the algorithm terminates, or some node prices will increase to infinity. This is the key idea that underlies the validity of the algorithm, and forms the basis for its proof of termination.

To describe formally the algorithm, consider the case where $P\ne (s)$ and $P$ has the form $P=(s,n_1,\ldots,n_k)$. We then denote by
$$\hbox{pred}(n_k)=n_{k-1}$$
the predecessor node of the terminal node $n_k$ in the path $P$. [In the case where $P= (s,n_1)$, we use the notation $\hbox{pred}(n_1)=s$.]  If the terminal node $n_k$ of $P$ is not deadend, we denote by
$\hbox{succ}(n_k)$
a downstream neighbor of $n_k$ that has minimal price:
$$\hbox{succ}(n_k)\in\arg\min_{\{j\,|\,(n_k,j)\in{\cal A}\}}p_j.$$
If multiple downstream neighbors of $n_k$ have minimal price, the algorithm designates arbitrarily one of these neighbors as $\hbox{succ}(n_k)$.

The algorithm also uses a positive scalar $\e$. The choice of $\e$ does not affect the path produced by the algorithm (so we could use $\e=1$ for example), but the choice of $\e$ will play an important role in the weighted path construction algorithm of the next section. 

The rules by which the path $P$ and the prices $p_i$ are updated at each iteration are as follows.

\xdef\figpathfind{\figr}\figrnum\show{myfigure}

\topinsert
\centerline{\hskip0pc\epsfxsize = 5.7in \epsfbox{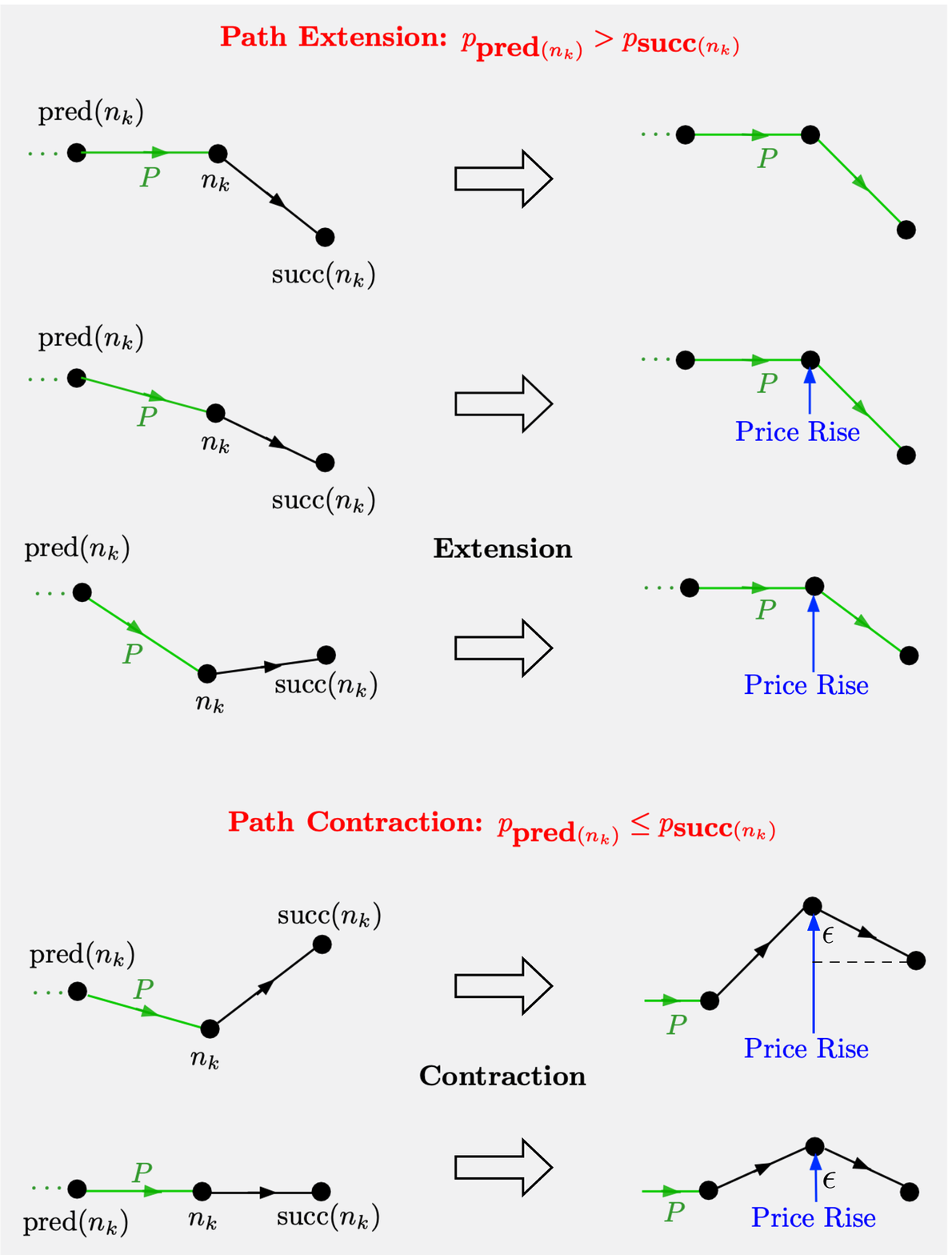}}%\vskip0pt
\fig{-0.5pc}{\figpathfind}{Illustration of the price levels of the terminal node $n_k$ of the path $P=(s,n_1,\ldots,n_k)$, and the price levels of its predecessor and its successor, before and after an extension or a contraction; cf.\ cases (c1) and (c2) of the algorithm description. In the case where
$p_{\hbox{pred}(n_k)}> p_{\hbox{succ}(n_k)},$
which corresponds to an extension, there is an increase of $p_{n_k}$ if the arc $\big(\hbox{pred}(n_k),n_k\big)$ is downhill. In the case where
$p_{\hbox{pred}(n_k)}\le  p_{\hbox{succ}(n_k)},$
which corresponds to a contraction, there is always an increase of $p_{n_k}$.}
\endinsert

\texshopbox{\pn {\bf Auction Algorithm Iteration for Unweighted Path Construction:} We distinguish three mutually exclusive cases.
\vskip0.5pc%}\texshopboxnt{
\pn{\bf (a)} {\it $P=(s)$\/}: We then set the price $p_s$  to $\max\{p_s,\,p_{\hbox{succ}(s)}+\e\}$, and extend $P$ to $\hbox{succ}(s)$.%}\texshopboxnt{
\pn{\bf (b)} {\it $P=(s,n_1,\ldots,n_k)$ and node $n_k$ is deadend\/}: We then  set the price $p_{n_k}$ to $\infty$ (or a very high number for practical purposes),  and contract $P$ to $n_{k-1}$.%}\texshopboxnt{
\pn{\bf (c)} {\it $P=(s,n_1,\ldots,n_k)$ and node $n_k$ is not deadend\/}. We then consider the following two cases.%}\texshopboxnt{
\nitem{(1)} $p_{\hbox{pred}(n_k)}> p_{\hbox{succ}(n_k)}.$
We then extend $P$ to $\hbox{succ}(n_k)$ and set $p_{n_k}$ to the price of $\hbox{pred}(n_k)$ [thus making the arc $\big(\hbox{pred}(n_k),n_k\big)$ level and the arc $\big(n_k,\hbox{succ}(n_k)\big)$ downhill].%}\texshopboxnt{
\nitem{(2)} $p_{\hbox{pred}(n_k)}\le p_{\hbox{succ}(n_k)}.$
We then contract $P$ to $\hbox{pred}(n_k)$ and raise the price of $n_k$ to the price of $\hbox{succ}(n_k)$ plus $\e$ [thus making the arc $\big(\hbox{pred}(n_k),n_k\big)$ uphill and the arc $\big(n_k,\hbox{succ}(n_k)\big)$ downhill].
}

The algorithm terminates once the destination becomes the terminal node of $P$. We will show that eventually the algorithm terminates, under our standing assumption that there is at least one path from the origin to the destination.

The contraction/extension mechanism of the algorithm may be interpreted as a competitive market process: we can view $\hbox{pred}(n_k)$ as being in competition with the downstream nodes of $n_k$ for becoming the next terminal node of path $P$, after $n_k$. In particular, the terminal node of $P$ moves to the node that offers minimal price [with ties that involve $\hbox{pred}(n_k)$ broken in favor of $\hbox{pred}(n_k)$ in order to maintain the downhill path property].

Figure \figpathfind\ illustrates the extension and contraction mechanism of case (c) above, and shows how the attendant price rises maintain the downhill path property of the algorithm throughout its operation. In particular, the initial path $P=(s)$ satisfies the downhill path property trivially, since it contains no arcs. Furthermore, using Fig.\ \figpathfind\ and the algorithm description, we can verify that if $P$ and the node prices satisfy  the downhill path property at the beginning of an iteration, then the new path and node prices at the beginning of the next iteration also  satisfy  the downhill path property. 

\vskip-0.5pc
 \subsection{Algorithm Justification}
\vskip-0.5pc

\pn We will prove that eventually the destination will become the terminal node of $P$, at which time the algorithm will terminate. To this end we argue by contradiction and we use our assumption that there is at least one path from the origin to the destination.

Indeed, suppose, to arrive at a contradiction, that the algorithm does not terminate. Then, since the path $P$ does not contain a cycle and hence cannot extend indefinitely, the algorithm must perform an infinite number of contractions. Let ${\cal N}_\infty$ be the nonempty set of nodes whose price increases infinitely often due to a contraction (and hence have a price that increases to $\infty$), and let $\ol{\cal N}_\infty=\{i\mid i\notin {\cal N}_\infty\}$ be the complementary set of nodes whose price increases due to a contraction finitely often (and hence do not become the terminal node of $P$ after some iteration). Clearly, by the rules of the algorithm, there is no arc connecting a node of  ${\cal N}_\infty$ to a node of $\ol{\cal N}_\infty$. Moreover, $t$ clearly belongs to $\ol{\cal N}_\infty$, and we  claim that $s$ belongs to ${\cal N}_\infty$. Indeed, if $s\in\ol{\cal N}_\infty$ there would exist a  subpath $P'=(s,n_1,\ldots,n_k)$  such that the nodes $s,n_1,\ldots,n_{k-1}$ belong to $\ol {\cal N}_\infty$, the last node $n_k$ belongs to ${\cal N}_\infty$, and $P'$ is the initial portion of $P$ for all iterations after finitely many. Since  $n_k$ will be the terminal node of $P$ infinitely often, it follows that  $n_{k-1}$ will be the predecessor $\hbox{pred}(n_k)$ of $n_k$ infinitely often, while the price of $n_k$ increases to infinity and the price of $n_{k-1}$ stays finite. By the rules of the algorithm, this is not possible. Thus we must have $s\in {\cal N}_\infty$, $t\in \ol{\cal N}_\infty$, and no arc connecting a node of  ${\cal N}_\infty$ to a node of $\ol{\cal N}_\infty$. This contradicts the assumption that  there is a path from $s$ to $t$, and shows that the algorithm will terminate.

We summarize the preceding arguments in the following proposition.

\xdef\propapc{\propn}\propnum\show{myproposition}

\texshopbox{\proposition{\propapc:} If there exists at least one path from
the origin to the destination, the APC algorithm terminates
with a path from $s$ to $t$.
Otherwise the algorithm never terminates and we have $p_i\to\infty$ for all nodes $i$ in a subset ${\cal N}_\infty$ that contains $s$.}

\vskip-1.5pc
\section{Path Planning with Arc Weights}
\vskip-0.5pc
\pn We will now introduce a generalization of our path planning algorithm, which we call {\it auction/weighted path construction} (AWPC for short). The algorithm  incorporates a length (or weight) $a_{ij}$ for every arc $(i,j)$, and aims to provide a path with near-minimum total length. Each length $a_{ij}$ encodes a measure of desirability of including arc $(i,j)$ into a path from the origin to the destination. The arc lengths serve to provide a bias towards producing paths with small total length. In fact in many cases (but not always) the algorithm produces shortest paths with respect to the given lengths. We require that {\it all cycles have nonnegative length\/}. By this we mean that for every cycle $(i,n_1,\ldots,n_k,i)$ we have 
$$a_{i,n_1}+a_{n_1n_2}+\cdots+a_{n_{k-1}n_k}+a_{n_ki}\ge 0.\eqnum\show{csao}$$
This is a common assumption in shortest path problems.\footnote{\dag}{\ninepoint The earlier auction/shortest path algorithm of [Ber91] requires that all cycle lengths be strictly positive rather than nonnegative, which is often a significant restriction. Another important difference, which affects computational efficiency, is that there is no $\e$ parameter in the algorithm of [Ber91]. Indeed, this algorithm is closely related to the so called ``naive auction algorithm," which is the auction algorithm for the assignment problem with $\e=0$.}

The AWPC algorithm maintains and updates a directed path $P=(s,n_1,\ldots,n_k)$ and a price $p_i$ for each node $i$. Extending the terminology of the preceding section, we say that under the current set of prices and lengths an arc $(i,j)$ is:
\nitem{(a)} {\it Downhill\/}: If $p_i>a_{ij}+p_j$.
\nitem{(b)} {\it Level\/}: If $p_i=a_{ij}+p_j$.
\nitem{(c)} {\it Uphill\/}: If $p_i<a_{ij}+p_j$.
\smskip
\pn As earlier, we denote by
$$\hbox{pred}(n_k)=n_{k-1}$$
the predecessor node of the terminal node $n_k$ in the path $P$. [In the case where $P= (s,n_1)$, we let $\hbox{pred}(n_1)=s$.]  If the terminal node $n_k$ of $P$ is not deadend, we denote by
$\hbox{succ}(n_k)$
a downstream neighbor $j$ of $n_k$ for which $a_{n_kj}+p_j$ is minimized:
$$\hbox{succ}(n_k)\in\arg\min_{\{j\,|\,(n_k,j)\in{\cal A}\}}\{a_{n_kj}+p_j\}.$$
If multiple downstream neighbors of $n_k$ attain the minimum, the algorithm designates arbitrarily one of these neighbors as $\hbox{succ}(n_k)$.

Note that when $a_{ij}=0$ for all arcs $(i,j)$, the preceding definitions coincide  with the ones given in the preceding section. Indeed when $a_{ij}=0$, the AWPC algorithm to be presented coincides with the APC algorithm of Section 2.

The AWPC algorithm maintains a directed path $P=(s,n_1,\ldots,n_k)$ that  starts at the origin and consists of distinct nodes. The path is either the degenerate path $P=(s)$, or it ends at some node $n_k\ne s$, which, as earlier,  is called the {\it terminal node} of $P$. Each iteration of the algorithm starts with a path and a price for each node, which are updated during the iteration. The algorithm starts with the degenerate path $P=(s)$, and the initial prices are  arbitrary.  It terminates when the destination  becomes the terminal node of $P$. 

We will now describe the rules by which the path and the prices are updated.
At any one iteration the algorithm starts with a path $P$ and a scalar price $p_i$ for each node $i$. At the end of the iteration a new path $\overline P$ is obtained from $P$ through a contraction or an extension as earlier. For iterations where the algorithm starts with the degenerate path $P=(s)$, only an extension is possible, i.e., $P=(s)$ is replaced by a path of the form $\overline P=(s,n_1)$.
Also the price of the terminal node of $P$ is increased just before a contraction, and in some cases, just before an extension. The amount of price rise is determined by a scalar parameter $\e>0$.

The algorithm terminates when the destination becomes the terminal node of $P$. The rules by which the path $P$ and the prices $p_i$ are updated at every iteration prior to termination are as follows.

\texshopboxnb{\pn {\bf Auction Algorithm Iteration for Weighted Path Construction:} We distinguish three mutually exclusive cases.
\vskip0.5pc
\pn{\bf (a)} {\it $P=(s)$\/}: We then set the price $p_s$  to $\max\{p_s,\,a_{s\hbox{succ}(s)}+p_{\hbox{succ}(s)}+\e\}$, and extend $P$ to $\hbox{succ}(s)$.%}\texshopboxnt{\pn
\pn{\bf (b)} {\it $P=(s,n_1,\ldots,n_k)$ and node $n_k$ is deadend\/}: We then  set the price $p_{n_k}$ to $\infty$ (or a very high number for practical purposes),  and contract $P$ to $n_{k-1}$. %}\texshopboxnt{\pn
\pn{\bf (c)} {\it $P=(s,n_1,\ldots,n_k)$ and node $n_k$ is not deadend\/}. We then consider the following two cases.%}\texshopboxnt{\pn
\nitem{(1)} $p_{\hbox{pred}(n_k)}> a_{\hbox{pred}(n_k)n_k}+a_{n_k\hbox{succ}(n_k)}+p_{\hbox{succ}(n_k)}.$
We then extend $P$ to $\hbox{succ}(n_k)$ and set $p_{n_k}$ to 
$$p_{\hbox{pred}(n_k)}-a_{\hbox{pred}(n_k)n_k},$$}\texshopboxnt{\pn 
[thus making the arc $\big(\hbox{pred}(n_k),n_k\big)$ level and the arc $\big(n_k,\hbox{succ}(n_k)\big)$ downhill].%}\texshopboxnt{\pn
\nitem{(2)} $p_{\hbox{pred}(n_k)}\le  a_{\hbox{pred}(n_k)n_k}+a_{n_k\hbox{succ}(n_k)}+p_{\hbox{succ}(n_k)}.$
We then contract $P$ to $\hbox{pred}(n_k)$ and raise the price of $n_k$ to 
$$a_{n_k\hbox{succ}(n_k)}+p_{\hbox{succ}(n_k)}+\e$$
[thus making the arc $\big(\hbox{pred}(n_k),n_k\big)$ uphill and the arc $\big(n_k,\hbox{succ}(n_k)\big)$ downhill].
}

A downhill/level/uphill type of interpretation, similar to Fig.\ \figpathfind, applies to this algorithm as well [the relative heights of the prices of nodes $\hbox{pred}(n_k)$, $\hbox{succ}(n_k)$, and $n_k$, indicated in Fig.\ \figpathfind\ should incorporate the arc lengths $a_{\hbox{pred}(n_k)n_k}$ and $a_{n_k\hbox{succ}(n_k)}$, as in the preceding algorithm description]. There is a price increase of $n_k$ in the case of a contraction, and also in the case of an extension if the arc $\big(\hbox{pred}(n_k),n_k\big)$ is downhill. However, the conditions for an arc $(i,j)$ to be downhill, level, or uphill involve the arc lengths $a_{ij}$.
Similar to our earlier arguments, it can be seen that $P$ and the prices $p_i$ satisfy the following downhill path property at the start of each iteration for which $P\ne (s)$. 

\texshopbox{\pn {\bf Downhill Path Property:}
\pn For the path $P=(s,n_1,\ldots,n_k)$ mantained by the algorithm, the last arc $(n_{k-1},n_k)$ of $P$ is level following a contraction to $n_k$ and is downhill following an extension to $n_k$. All other arcs of $P$ are level.}

\xdef\figureawsp{\figr}\figrnum\show{myfigure}

A consequence of this property (and our assumption that all cycles have nonnegative length) is that when an extension occurs, a cycle cannot be created, in the sense  that the terminal node $n_k$ is different than all the predecessor nodes $s,n_1,\ldots,n_{k-1}$ on the path $P=(s,n_1,\ldots,n_k)$.
Thus, assuming that there is at least one path from the origin to the destination, it can be shown that eventually the destination will become the terminal node of $P$, at which time the algorithm will terminate. The proof is essentially identical to the proof we gave earlier for the case of zero lengths in Prop.\ \propapc. Figure \figureawsp\ illustrates the
algorithm .

\midinsert
\centerline{\hskip0pc\epsfxsize = 5.2in \epsfbox{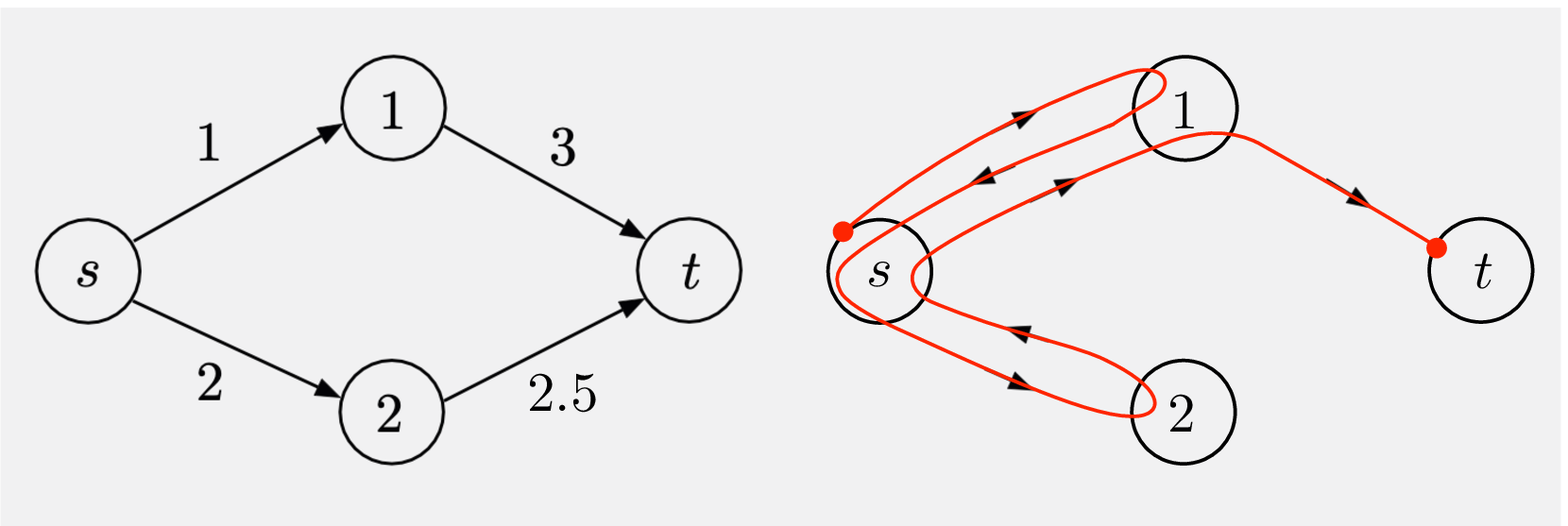}}
\vskip0pt
\eightpoint
\def\tablerule{\noalign{\hrule}}
$$\vbox{\offinterlineskip
\hrule
\halign{\vrule\hfill \ #\ \hfill &\vrule\hfill \ #\ \hfill 
&\vrule\hfill \ #\ \hfill &\vrule\hfill 
\ #\ \hfill \vrule\cr
&&\cr
{\bf Iteration \#\lower2ex\hbox{\ }\raise4ex\hbox{\ }}&\hbox{\bf
Path $P$ prior}&\hbox{\bf Price vector $(p_s,p_1,p_2,p_t)$} &\hbox{\bf Type of iteration}\cr
{\bf \lower2ex\hbox{\ }}&\hbox{\bf
to iteration}&\hbox{\bf  prior to iteration} &\cr
\tablerule\cr &&\cr
{1\lower2ex\hbox{\ 
}\raise4ex\hbox{\ }\hbox{\ }}&$(s)$&$(\underline{0},0,0,0)$&Extension to 1\cr 
{2\lower2ex\hbox{\
}\raise2ex\hbox{\ }\hbox{\ }}&$(s,1)$&$(1+\e,\underline{0},0,0)$&Contraction to $s$\cr 
{3\lower2ex\hbox{\
}\raise2ex\hbox{\ }\hbox{\ }}&$(s)$&$(\underline{1+\e},3+\e,0,0)$&Extension to 2\cr 
{4\lower2ex\hbox{\
}\raise2ex\hbox{\ }\hbox{\ }}&$(s,2)$&$(2+\e,3+\e,\underline{0},0)$&Contraction to $s$\cr 
{5\lower2ex\hbox{\
}\raise2ex\hbox{\ }\hbox{\ }}&$(s)$&$(\underline{2+\e},3+\e,2.5+\e,0)$&Extension to 1\cr
{6\lower2ex\hbox{\
}\raise2ex\hbox{\ }\hbox{\ }}&$(s,1)$&$(4+2\e,\underline{3+\e},2.5+\e,0)$&Extension to $t$\cr
{7\lower2ex\hbox{\
}\raise2ex\hbox{\ }\hbox{\ }}&$(s,1,t)$&$(4+2\e,3+2\e,2.5+\e,\underline{0})$&Termination\cr
&&\cr \tablerule\cr}}$$
\fig{-1pc}{\figureawsp:}{A four-node example, with the arc lengths  shown next to the arcs in the left-side figure. The right-side figure and the table trace the steps of  the AWPC algorithm starting with $P=(s)$ and
all initial prices equal to 0. The price of the terminal node of the path is underlined. The trajectory of the AWPC algorithm shown in the table corresponds to values of $\e\le3$. The final path obtained is the shortest path $(s,1,t)$. If instead $\e> 3$, the algorithm will still find the shortest path and faster: it will first perform an extension to 1, setting $p_s=1+\e$. It will then perform an extension to $t$, since the condition
$$1+\e=p_{\hbox{pred}(1)}> a_{\hbox{pred}(1)1}+a_{1\hbox{succ}(1)}+p_{\hbox{succ}(1)}=1+3+0=4$$
is satisfied, and terminate. The final path will be $P=(s,1,t)$ and the final price vector will be 
$$(p_s,p_1,p_2,p_t)=(1+\e,1+\e,0,0),$$
with the arc $(s,,1)$ being balanced and the arc $(1,t)$ being downhill.
Note also that generally, there is no guarantee that the AWPC algorithm will find a shortest path for all initial prices and values of $\e$. However, it will find a shortest path if $\e$ is sufficiently small, provided the initial prices satisfy a certain $\e$-complementary slackness constraint that will be given in Section 3.2.}\endinsert

\subsection{The Role of the Parameter $\e$ - Convergence Rate and Solution Accuracy Tradeoff, and $\e$-Scaling}
\vskip-0.5pc

\pn In auction algorithms, it is common to use a positive $\e$ parameter to regulate the size of price rises. In the AWPC algorithm, $\e$ is used to provide an important tradeoff between the ability of the algorithm to construct paths with near-minimum  length, and its rate of convergence. Generally, as $\e$ becomes smaller the quality of the path produced improves, as we will show with examples and analysis in what follows. On the other hand a small value of $\e$ tends to slow down the algorithm. 

The tradeoff between speed of convergence and accuracy of solution that is embodied in the choice of $\e$ was recognized in the original proposal of the auction algorithm for the assignment problem [Ber79], and the approach of {\it $\e$-scaling} was proposed to deal with it. In this approach we start the auction algorithm with a relatively large value of $\e$, to obtain quickly rough estimates for appropriate values of the node prices, and then we progressively reduce $\e$ to refine the node prices and eventually obtain an optimal solution. The use of $\e$-scaling also allows the option of stopping the algorithm, with a less refined solution, if the allotted time for computation is limited.

In the context of the AWPC algorithm, $\e$-scaling is implemented by running the algorithm with a relatively large value of $\e$ to estimate ``good" prices, at least for a subset of ``promising" transit nodes from $s$ to $t$, and then progressively refining the assessment of the ``promise" of these nodes. This is done by rerunning the algorithm with smaller values of $\e$, while using as initial prices at each run the final prices of the previous run. It is well-known that $\e$-scaling improves the computational complexity of auction algorithms,\footnote{\dag}{\ninepoint See the papers [BeE88], [Ber88], the book [Ber98], and the references quoted there, for polynomial complexity analyses of auction algorithms for the assignment problem and other related problems.} and it can be similarly applied to the AWPC algorithm, as we will discuss shortly.

In what follows in this section, we will use examples to illustrate how the choice of $\e$ affects the rate of convergence of the AWPC algorithm, as well as the error from optimality of the path that it produces.

\xdef\examplecycle{\exampl}\examplnum\show{myexample}

\xdef\figexamplecycle{\figr}\figrnum\show{myfigure}

\topinsert
\centerline{\hskip0pc\epsfxsize = 4.2in \epsfbox{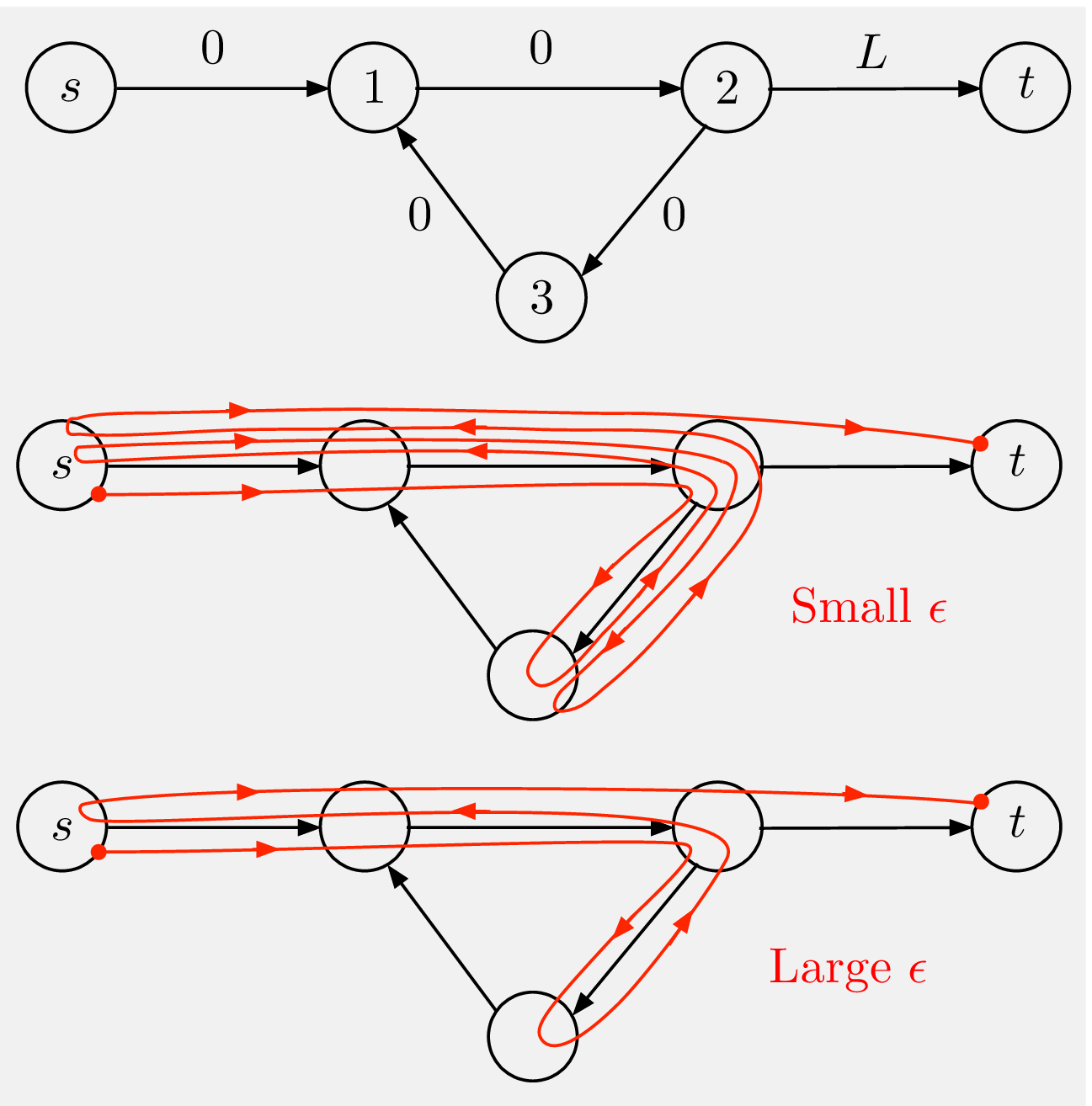}}%\vskip0pt
\vskip-8pt
\eightpoint
\def\tablerule{\noalign{\hrule}}
$$\vbox{\offinterlineskip
\hrule
\halign{\vrule\hfill \ #\ \hfill &\vrule\hfill \ #\ \hfill 
&\vrule\hfill \ #\ \hfill &\vrule\hfill 
\ #\ \hfill \vrule\cr
&&\cr
{\bf Iteration \#\lower2ex\hbox{\ }\raise4ex\hbox{\ }}&\hbox{\bf
Path $P$ prior}&\hbox{\bf Price vector $(p_s,p_1,p_2,p_3,p_t)$} &\hbox{\bf Type of iteration}\cr
{\bf \lower2ex\hbox{\ }}&\hbox{\bf
to iteration}&\hbox{\bf  prior to iteration} &\cr
\tablerule\cr &&\cr
{1\lower2ex\hbox{\ 
}\raise4ex\hbox{\ }\hbox{\ }}&$(s)$&$(\underline{0},0,0,0,0)$&Extension to 1\cr 
{2\lower2ex\hbox{\
}\raise2ex\hbox{\ }\hbox{\ }}&$(s,1)$&$(\e,\underline{0},0,0,0)$&Extension to 2\cr 
{3\lower2ex\hbox{\
}\raise2ex\hbox{\ }\hbox{\ }}&$(s,1,2)$&$(\e,\e,\underline{0},0,0)$&Extension to 3\cr
{4\lower2ex\hbox{\
}\raise2ex\hbox{\ }\hbox{\ }}&$(s,1,2,3)$&$(\e,\e,\e,\underline{0},0)$&Contraction to 2\cr
{5\lower2ex\hbox{\
}\raise2ex\hbox{\ }\hbox{\ }}&$(s,1,2)$&$(\e,\e,\underline{\e},2\e,0)$&Contraction to 1\cr
{6\lower2ex\hbox{\
}\raise2ex\hbox{\ }\hbox{\ }}&$(s,1)$&$(\e,\underline{\e},3\e,2\e,0)$&Contraction to $s$\cr
{7\lower2ex\hbox{\
}\raise2ex\hbox{\ }\hbox{\ }}&$(s)$&$(\underline{\e},4\e,3\e,2\e,0)$&Extension to 1\cr
{8\lower2ex\hbox{\
}\raise2ex\hbox{\ }\hbox{\ }}&$(s,1)$&$(5\e,\underline{4\e},3\e,2\e,0)$&Extension to 2\cr
{9\lower2ex\hbox{\
}\raise2ex\hbox{\ }\hbox{\ }}&$(s,1,2)$&$(5\e,5\e,\underline{3\e},2\e,0)$&Extension to $t$ if $2\e>L$\cr
{\lower2ex\hbox{\
}\raise2ex\hbox{\ }\hbox{\ }}&&&Extension to $2$ otherwise\cr
{$\ldots$\lower2ex\hbox{\
}\raise2ex\hbox{\ }\hbox{\ }}&$\ldots$&$\ldots$&Continue until $p_3>L$\cr
\tablerule\cr}}$$
\fig{-1.5pc}{\figexamplecycle}{The shortest path problem of Example \examplecycle\ (top part of the figure). The arc lengths are shown next to the arcs [all lengths are equal to 0, except for the length of arc $(2,t)$ which has a large length $L$]. There is only one point where the algorithm can go wrong, at node 2 where there is a choice between going to $t$ or going to 3. The only $s$-to-$t$ path is $(s,1,2,t)$, but if $\e$ is very small, the algorithm explores the possibility of reaching the destination through node 3 for many iterations, while repeating the cycle
$s\to1\to2\to3\to2\to1\to s\to1\ldots$
(middle part of the figure). On the other hand, if $2\e>L$, then at iteration 9, following an extension to node 2, the successor to node 2 is $t$, the algorithm compares the prices of nodes 3 and $t$, performs an extension to $t$, and terminates (bottom part of the figure).
}
\endinsert

\beginexample{\examplecycle\ ($\e$-Scaling and the Effect of Small Cycles)}\pn This is an example of a shortest path problem where there is a cycle of relatively small length. It involves that graph of the top figure in Fig.\ \figexamplecycle. The cycle consists of nodes 1, 2, and 3, and has length 0 [the algorithm's behavior is similar when the cycle has positive length that is small relative to $L$, the length of the last arc of the unique $s$-to-$t$ path]. Such cycles slow down the algorithm, when $\e$ has a small value. Indeed, it can be seen from the table of   Fig.\ \figexamplecycle\ that for small values of $\e$ and initial prices equal to 0, the algorithm repeats the cycle  
$$s\to1\to2\to3\to2\to1\to s\to1\cdots$$
until the prices of nodes 2 and 3 reach levels $p_2> L$ and $p_3>L$, so that the arc $(2,t)$ becomes downhill and an extension from 2 to $t$ is performed. The number of cycles for this to happen depends on $\e$ and is roughly proportional to $L/\e$, so for small values of $\e$, the computation is nonpolynomial (see the middle part of the figure). On the other hand, it can be seen from the table of Fig.\ \figexamplecycle\ that when $\e$ is large enough so that $3+4\e>L$, the algorithm moves to $t$ once it reaches node 2 for the second time, after 8 iterations.

It can be shown that with an $\e$-scaling scheme, whereby $\e$ is reduced by a certain factor between successive runs of the algorithm, the computation is polynomial, proportional to $\log L$, rather than $L$. This is a common property of auction algorithms; see the book [Ber98] and the references quoted there.  Later, in Section 3.3, we will discuss the use of $\e$-scaling and its use to provide convergence acceleration as well as exact shortest path solutions, albeit with an $\e$-complementary restriction in the choice of the initial prices.
\endexample

In the preceding example, the poor performance of the algorithm is caused by the presence of a cycle with small length. The next example illustrates how a similar phenomenon can also occur in acyclic graphs involving many-node paths.

\xdef\examplechain{\exampl}\examplnum\show{myexample}
\xdef\figexamplechain{\figr}\figrnum\show{myfigure}

\beginexample{\examplechain\ (Convergence Rate and Solution Accuracy Tradeoff)}\pn Consider a graph involving a long chain of nodes that starts at the origin and ends at the destination, as shown in Fig.\ \figexamplechain. We assume that the initial prices are all equal to 0. Then it can be verified that for large values of $\e$, the algorithm will terminate with the suboptimal path $(s,1,2,\ldots,n,t)$; in fact for $\e>n$, it will terminate in $n+1$ iterations through the sequence of extensions 
$s\to1\to2\to\cdots \to n\to t.$
 It can also be verified, by tracing the steps of the algorithm that for small values of $\e$, the algorithm will find the optimal path $(s,2',t)$, but will need a large number of iterations (proportional to $n^2$) to do so. A suitable $\e$-scaling scheme can find the optimal path in $O(n\log n)$ iterations.
\endexample 

\topinsert
\centerline{\hskip0pc\epsfxsize = 5.2in \epsfbox{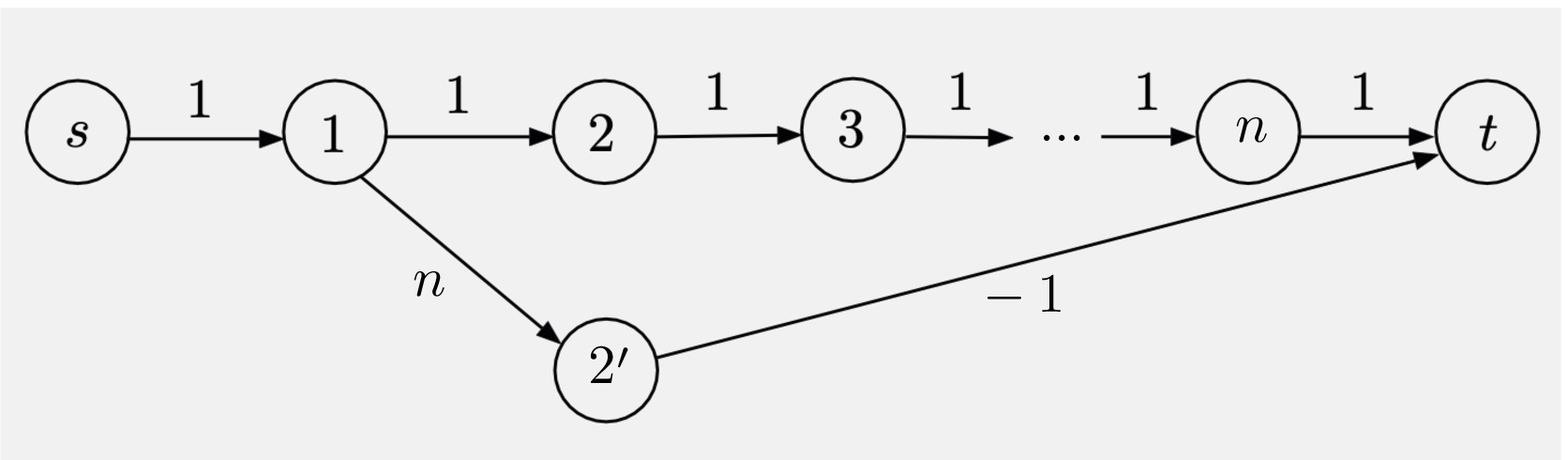}}%\vskip0pt
\vskip0pt
\fig{-1pc}{\figexamplechain}{The graph of the shortest path problem of Example \examplechain. The arc lengths are shown next to the arcs. All lengths are equal to 1, except for  the length of arc $(2',t)$ which is equal to  $-1$ and the length of arc $(1,2')$ which is equal to  $n$, the number of intermediate nodes in the top path from $s$ to $t$. The shortest path is $(s,2',t)$, but for large values of $\e$, the algorithm will terminate with the suboptimal path $(s,1,2,\ldots,n,t)$.}
\endinsert
\vskip-0.5pc

\subsection{Initial Price Selection and Shortest Distances}

\pn An important issue in the AWPC algorithm is the choice of the initial prices. We will argue that the algorithm operates effectively, in the sense that it terminates fast and with small error from path optimality, if the initial prices are close to the shortest distances under the given set of lengths. Indeed, the lengths $\{a_{ij}\}$ define the shortest distances, denoted by  $d^*_i$, from the nodes $i$ to the destination $t$. These shortest distances satisfy $d^*_t=0$ and for all $i\ne t$,
$$d^*_i=\min_{\{j\mid (i,j)\in{\cal A}\}}\{a_{ij}+d^*_j\};$$
this is an instance of the fundamental dynamic programming/Bellman equation. Suppose that we choose for all $i$ a  price $p_i$ that is exactly equal to $d^*_i$. Then it can be verified that starting from an arbitrary origin $i$, the algorithm generates a shortest path from $i$ to $t$ through a sequence of extensions over level arcs, without any intervening contractions. The price differential $p_s-p_t$ is equal to the total length of the path produced by the algorithm, which is shortest.

If the initial prices $p_i$ are not equal to the shortest distances $d^*_i$, the price differential $p_s-p_t$ provides an upper bound to the total length of the final path $P$ produced by the algorithm:
$$ L_P\le p_s-p_t,\xdef\pathineq{\lab}\eqnum\show{csao}$$
To see this, note that we have
$a_{sn_1}=p_s-p_{n_1}$, $a_{n_{i-1}n_i}=p_{n_{i-1}}-p_{n_i}$ for all $i=1,\ldots,k$,  and $a_{n_kt}\le p_{n_k}-p_t$, since in view of the downhill property, all the arcs of $P$ are level except possibly for the last one, which is downhill. It follows that
$$L_P=a_{sn_1}+a_{n_1n_2}+\cdots+a_{n_{k-1}n_k}+a_{n_kt}\le(p_s-p_{n_1})+(p_{n_1}-p_{n_2})+\cdots+(p_{n_{k-1}}-p_{n_k})+(p_{n_k}-p_t)=p_s-p_t,$$
thus verifying Eq.\ \pathineq. 

 If all the arcs that  do not belong to $P$ are level or uphill upon termination, then it can be seen that $P$  has minimum total length. More generally, given any path $\hat P$ from $s$ to $t$, the total lengths $L_P$ and $L_{\hat P}$ of $P$ and $\hat P$ satisfy upon termination
$$L_P\le L_{\hat P}+\sum_{(i,j)\in \hat P}r_{ij},\xdef\suboptineq{\lab}\eqnum\show{csao}$$
where for any arc $(i,j)$ we denote by $r_{ij}$
the amount by which $a_{ij}$ must be increased to make the arc $(i,j)$ level if it is downhill:
$$r_{ij}=\max\{0,\,p_i-a_{ij}-p_j\}.\xdef\discrepancy{\lab}\eqnum\show{csao}$$
To see this note that for any path $\hat P$ that starts at $s$ and ends at $t$, we have, using Eq.\ \pathineq,
$$L_{\hat P}=\sum_{(i,j)\in \hat P}a_{ij}=p_s-p_t+\sum_{(i,j)\in \hat P}(a_{ij}+p_j-p_i)\ge L_P+\sum_{(i,j)\in \hat P}(a_{ij}+p_j-p_i)\ge L_P-\sum_{(i,j)\in \hat P}r_{ij},$$
thus verifying Eq.\ \suboptineq. 

Note that the preceding calculations are valid and the bound \suboptineq\ holds even if multiple arcs of $P$ are downhill (as long as there are no uphill arcs), so that we continue to have $p_s-p_t\ge L_{P}$ [cf.\ Eq.\ \pathineq]. This observation is relevant for variants of the algorithm for which different forms of the downhill path property hold; see e.g., Section 3.3.

We call the scalar $r_{ij}$ the {\it discrepancy} of arc $(i,j)$. The discrepancies provide an upper bound to the degree of suboptimality of the path obtained upon termination, as per Eq.\ \suboptineq. In particular, if $b$ is an upper bound to the arc discrepancies upon termination of the algorithm, then Eq.\ \suboptineq\ implies that
$$L_P\le L_{\hat P}+(n+1)b,\xdef\suboptineqtwo{\lab}\eqnum\show{csao}$$
where $n$ is the number of nodes other than $s$ and $t$. The reason is that a path from $s$ to $t$ can contain at most $(n+1)$ arcs, each having a discrepancy that is at most $b$. \old{Note that a discrepancy of an arc that is outside the path $P$ is created when there is a contraction (in which case the discrepancy is equal to $\e$), but it can also be created when there is a price rise through an extension [in which case arcs $(n_k,j)$ with $j\ne \hbox{succ}(n_k)$ may become downhill].}

An interesting empirical observation is that when the algorithm creates a new downhill arc $(i,j)$ that lies outside $P$, the corresponding discrepancy $r_{ij}$ becomes equal to $\e$ or a small multiple of $\e$. A reasonable conjecture is that if all the discrepancies $r_{ij}$ are initially bounded by a small multiple of $\e$, then the                                                                                                                                                                                       path produced by the algorithm upon termination is shortest to within a small multiple of $n\e$, where $n$ is the number of nodes. This bears similarity to auction algorithms for assignment and network flow problems, where the solution obtained can be proved to be optimal to within $n\e$. 

\subsection{Using $\e$-Scaling to Find a Shortest Path - $\e$-Complementary Slackness}
\vskip-0.5pc

\pn As noted earlier, $\e$-scaling consists of starting the auction algorithm with a relatively large value of $\e$, to obtain rough estimates for appropriate values of the node prices, and then progressively reducing $\e$ to refine the node prices and eventually obtain an optimal solution.  The idea is to improve the solution accuracy, without incurring the nonpolynomial convergence rate penalty illustrated in Example \examplecycle. One possible approach is to simply use the final prices obtained by the AWPC algorithm for a given value of $\e$ as initial prices for running AWPC for a smaller value of $\e$. This is possible because the AWPC algorithm can use arbitrary initial prices. Empirically, this scheme seems to work well, but it does not offer a guarantee that it will yield a shortest path. As an example in the problem of Fig.\ \figexamplechain, if all prices are chosen to be 0 except for the price of node 2', which is chosen to be a high positive number, the AWPC algorithm will not find the shortest path for any value of $\e$.

We will now consider $\e$-scaling schemes that offer guarantees of obtaining a shortest path. These schemes modify the prices produced by the AWPC algorithm for a given value of $\e$, before applying the algorithm with a smaller value of $\e$.  To this end, we will require that the algorithm satisfies initially and maintains the following {\it $\e$-complementary slackness} condition ($\e$-CS for short).

\texshopbox{\pn {\bf $\e$-Complementary Slackness:}
\pn For a given $\e>0$, the prices $\{p_i\mid i\in{\cal N}\}$ and the path $P$ satisfy
$$p_i\le a_{ij}+p_j+\e,\qquad \hbox{for all arcs }(i,j),$$
i.e., every arc is uphill, or level, or  downhill by at most $\e$, and %}\texshopboxnt{\pn 
$$p_i\ge a_{ij}+p_j,\qquad \hbox{for all arcs $(i,j)$ of the path $P$},$$
i.e., every arc of $P$ is level or downhill (by at most $\e$).}

The notion of $\e$-CS is fundamental in the context of auction algorithms, and represents a relaxation of the classical complementary slackness condition of linear programming (see, e.g., Bertsimas and Tsitsiklis [BeT97]). Note that when $\e$-CS holds, the discrepancies $r_{ij}$ of Eq.\ \discrepancy\ are at most equal to $\e$, so {\it if the AWPC algorithm maintains $\e$-CS throughout its operation, it produces a path that is suboptimal by at most $(n+1)\e$} in view of Eq.\ \suboptineqtwo, and hence also optimal for $\e$ sufficiently small [$(n+1)\e$ should be less than the difference between the 2nd shortest path distance and the shortest path distance]. It is possible to modify the price rise rule when an extension occurs in the AWPC algorithm, so that it maintains $\e$-CS throughout its operation if started with prices satisfying $\e$-CS.\footnote{\dag}{\ninepoint This requires a simple modification of the amount of price increase in case (c1), where 
$$p_{\hbox{pred}(n_k)}> a_{\hbox{pred}(n_k)n_k}+a_{n_k\hbox{succ}(n_k)}+p_{\hbox{succ}(n_k)}.$$ 
We raise the price $p_{n_k}$ to the largest value
that satisfies $\e$-CS, while extending $P$ to $\hbox{succ}(n_k)$, rather than raising $p_{n_k}$ to the value $$p_{\hbox{pred}(n_k)}- a_{\hbox{pred}(n_k)n_k}.$$
Then the downhill path property holds in a modified form, whereby all arcs of $P$ are either level or downhill, with the last arc of $P$ being downhill following an extension, but still guarantees that no cycle is created via an extension.} On the other hand finding initial prices that satisfy $\e$-CS may be difficult, although it is possible when $a_{ij}\ge0$ for all arcs $(i,j)$, in which case we can take $p_i=0$ for all nodes $i$. Moreover, operating the algorithm to satisfy $\e$-CS is complicated when arc lengths change over time and on-line replanning is necessary.

Given the final set of prices and path obtained by the AWPC algorithm for a given value of $\e$, there is an important issue in $\e$-scaling: how to modify the prices of some of the nodes so that the resulting prices together with the degenerate path $(s)$ satisfy $\e'$-CS for a smaller positive value $\e'<\e$. Moreover the price modifications should be small in order for the new prices to be good starting points for rerunning the AWPC algorithm with the new value $\e'$.

There are algorithms for computing price modifications to satisfy $\e'$-CS for a smaller value $\e'<\e$ together with the degenerate path $(s)$, which will be discussed in future reports. Moreover, often such algorithms can take advantage of special structure of the problem's graph. This is true for example  in assignment problems, where the bipartite character of the graph allows great flexibility in the choice of the initial prices. In what follows in this section, we discuss the case of shortest path problems involving an acyclic graph, which arise prominently in on-line and off-line multistep lookahead minimization and tree search for reinforcement learning problems (see Section 3.4). 

\subsubsection{$\e$-Scaling in Acyclic Graphs}

\xdef\figacyclicgraph{\figr}\figrnum\show{myfigure}

\topinsert
\centerline{\hskip0pc\epsfxsize = 3.0in \epsfbox{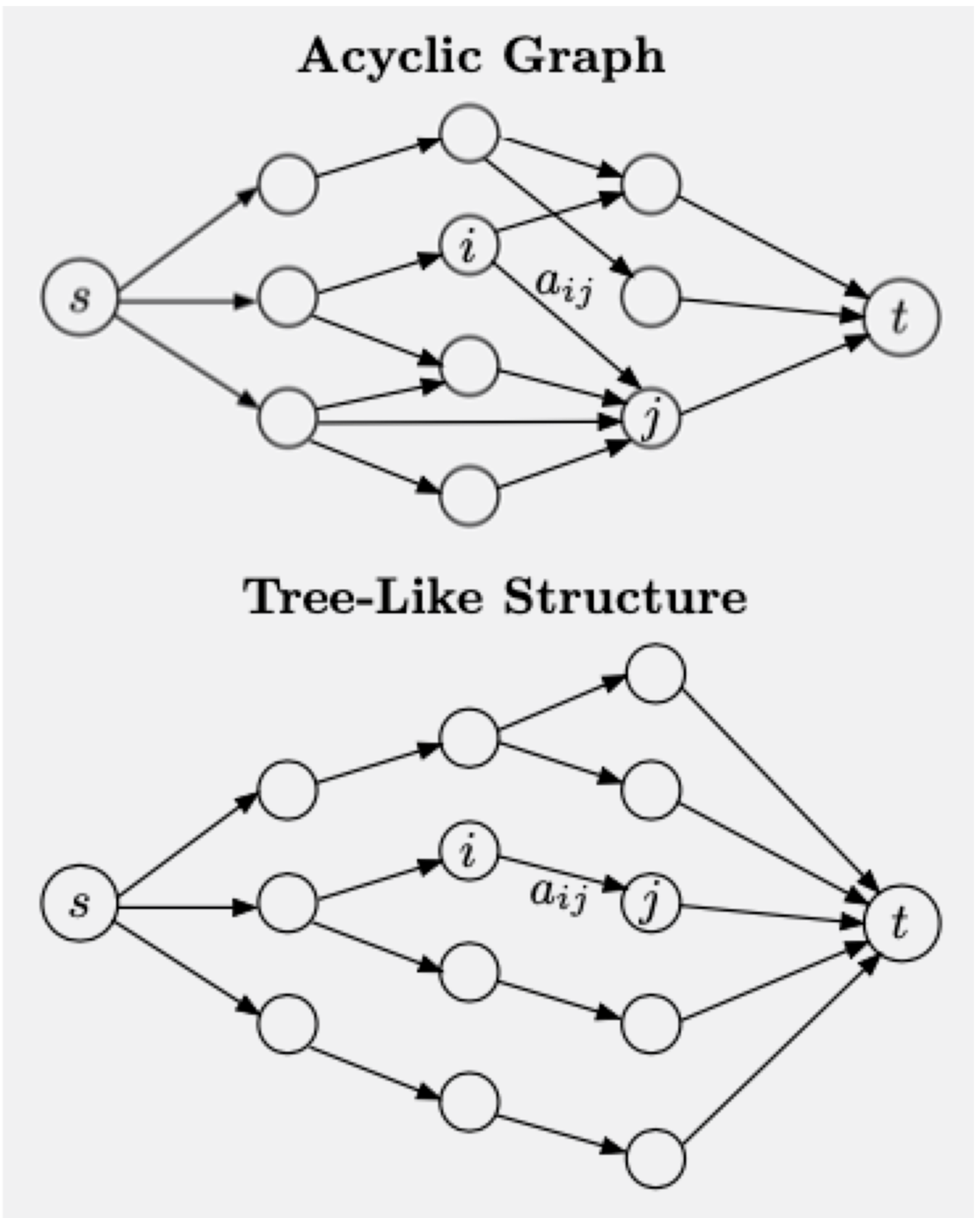}}%\vskip0pt
\fig{-0.5pc}{\figacyclicgraph}{Illustration of an acyclic graph involving paths that start at $s$ and end at $t$ and arc lengths $a_{ij}$. A  case of special interest in reinforcement learning is a tree-like structure, illustrated in the bottom figure, where nodes are grouped in layers, with arcs starting from one layer and ending at a node of the next layer, and  there is a single incoming arc to each node except $s$ and $t$.}
\endinsert

\pn Let us consider the case of an acyclic graph; cf.\ Fig.\ \figacyclicgraph. We are given arc lengths $a_{ij}$ and a set of prices $\{p_i\mid i\in {\cal N}\}$ that for a given positive $\e$ satisfy 
$$p_i\le a_{ij}+p_j+\e,\qquad \hbox{for all arcs }(i,j),\xdef\ecs{\lab}\eqnum\show{csao}$$
possibly resulting from application of the AWPC algorithm to the corresponding shortest path problem.
We want to find a set of prices $\{p_i'\mid i\in {\cal N}\}$ that for a given positive $\e'<\e$, satisfy 
$$p_i'\le a_{ij}+p_j'+\e',\qquad \hbox{for all arcs }(i,j),\xdef\ecsprime{\lab}\eqnum\show{csao}$$
satisfy $\e'$-CS together with the degenerate path $P=(s)$.
Thus Eq.\ \ecs\ requires that arcs, when downhill, are downhill by at most $\e$, while Eq.\ \ecsprime\ requires them to be downhill by at most $\e'$.

\xdef\figepsilonscaling{\figr}\figrnum\show{myfigure}

The idea is to start at $t$ and sequentially proceed backwards towards $s$, by delineating arcs $(i,j)$ that violate the condition \ecsprime\ and raising the price of $j$ and possibly the prices of some descendants of $j$ [since increasing $p_j$ may violate $\e'$-CS for nodes that lie downstream of $j$]. Thus we must check descendants of $j$ all the way to the destination $t$, and raise their prices  by whatever amounts are necessary to enforce the $\e'$-CS condition \ecsprime\ on arc $(i,j)$. Figure \figepsilonscaling\ provides an example. 

\topinsert
\centerline{\hskip0pc\epsfxsize = 3.8in \epsfbox{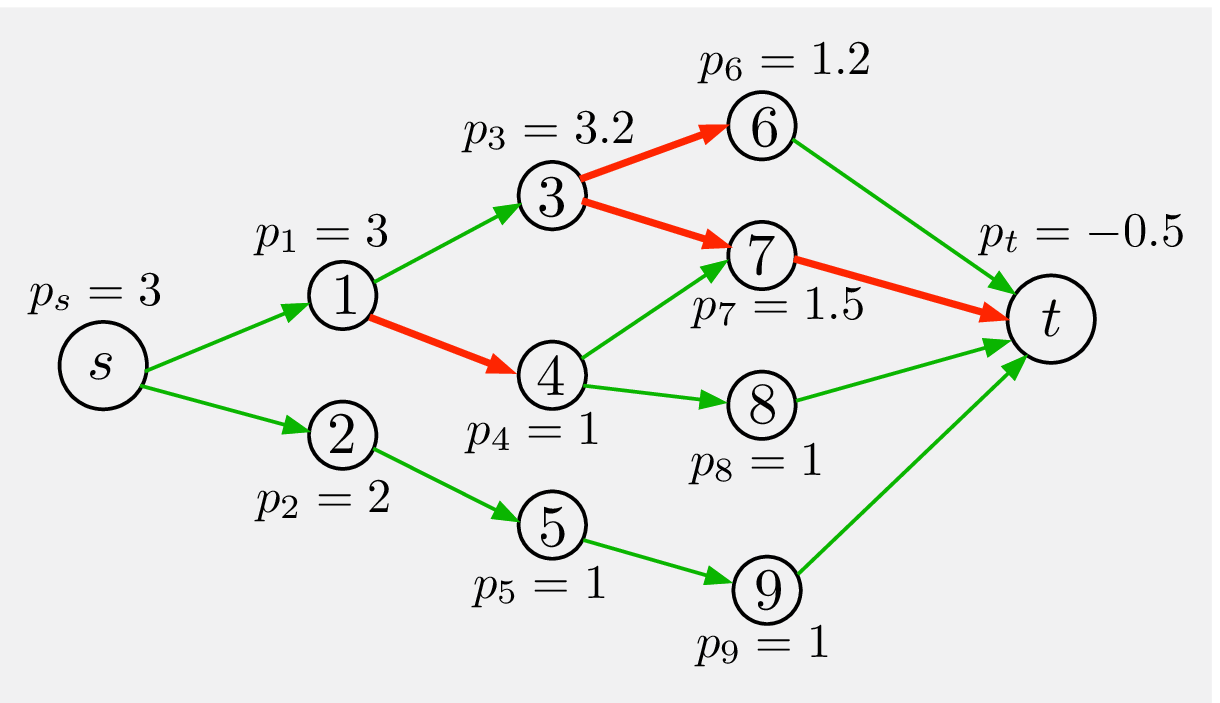}}
%\vskip0pt
\fig{-0.5pc}{\figepsilonscaling}{Illustration of $\e$-scaling in an acyclic graph, and the modifications needed to pass from the path $P=(s)$ and prices satisfying $\e$-CS to prices satisfying $\e'$-CS, with $\e'<\e$. All arcs have length equal to 1, and the prices $p_i$ are shown next to the nodes $i$. Let $\e=1$ and $\e'=1/2$. All arcs satisfy the $1$-CS condition \ecs, but arcs $(1,4)$, $(3,6)$, $(3,7)$, and $(7,t)$ (shown in red) violate the $0.5$-CS condition \ecsprime. We obtain prices satisfying  $0.5$-CS, by increasing the prices of the end nodes $t$, 6, 7, and 1 (in that order), and possibly their descendants in four iterations:
\nitem{(1)} Set $p_t=-0.5\uparrow 0$.
\nitem{(2)} Set $p_6=1.2\uparrow 1.7$ and $p_t=0\uparrow 0.2$.
\nitem{(3)} Set $p_7=1.5\uparrow 1.7$ (no need to increase $p_t$ further).
\nitem{(4)} Set $p_4=1\uparrow 1.5$ (no need to increase $p_7$ or $p_8$, and hence also $p_t$).
}
\endinsert

\subsection{Connections with Reinforcement Learning Methods}
\vskip-0.5pc
\pn The preceding analysis suggests that the initial prices $p_i$ should be chosen to be close to the shortest distances $d^*_i$, or more accurately, they should be chosen in a way that keeps the arcs nearly level or uphill, and minimizes the discrepancies given by Eq.\ \discrepancy. Of course we do not know the exact values $d^*_i$, but in a given application we may be able to use as initial prices approximate values, which may be obtained through a computationally inexpensive heuristic or other machine learning methods. 

In one possible off-line approach one may use data to train a neural network or other approximation architecture to learn approximations to the shortest distances $d_i^*$. The data may be obtained by using a shortest path algorithm and arc lengths that are similar to the ones of the given problem. The training should also aim to produce prices for which the discrepancies $r_{ij}$ given by Eq.\ \discrepancy\ are small. This objective can and should be encoded into the training problem.  It is also possible to train multiple neural networks to use for different patterns of arc lengths.

Another interesting context arises in a time-varying environment where some of the lengths may be changing as some arcs may become unavailable and new arcs may become available, while new instances of shortest path problems arise. In this case, an interesting possibility may be to update the initial prices using machine learning methodology, and a combination of off-line and on-line training with data. 

In this regard, we should mention that the use of reinforcement learning (RL) methods in conjunction with our path construction algorithms is facilitated by the fact that the initial prices are unrestricted. This makes our algorithms well-suited for large-scale and time-varying environments, such as data mining and transportation, where requests for solution of path construction problems arise continuously over time. Addressing the special implementation and machine learning issues in the context of such environments is an interesting subject for further research. 

There are also possibilities for incorporation of the AWPC algorithm within the RL methodology. Indeed, several RL methods rely on the computation of (nearly) shortest paths. Notable examples include multistep lookahead minimization and tree search methods. These methods involve evaluation of various decisions that are available at a given state by means of (nearly) shortest path calculations over an acyclic  decision tree, such as the one shown in the bottom part of Fig.\ \figacyclicgraph\ (see RL books such as the author's [Ber19], [Ber20a], [Ber22b], as well as the books by Powell [Pow11], and Sutton and Barto [SuB18]). The techniques of real-time dynamic programming, described in the papers by Korf [Kor90], and Barto, Bradtke, and Singh [BBS95], among many others, are also relevant in this context. 

The AWPC algorithm is inherently deterministic, but it can also be applied in stochastic multistep lookahead  contexts, where the popular simulation-based Monte Carlo tree search methods have been used widely. This can be done by replacing all steps of a multistep lookahead {\it except for the first} by deterministic approximations through the use of certainty equivalence (see e.g., the book [Ber19]).  The deterministic shortest path optimizations following the first step of lookahead involve an acyclic graph, and can be handled with the AWPC algorithm.  There is good reason for taking into account the stochastic nature of the first step, in order to maintain the connection of the lookahead minimization with Newton's method for solving the Bellman equation, as has been explained in the book [Ber22b].

We also note that heuristics, including  shortest path-type methods like $A^*$, arise prominently within rollout algorithms. This is a popular class of RL methods that has received a lot of attention as an effective and easily implementable (suboptimal) methodology for discrete and combinatorial optimization problems; see the books [Ber19], [Ber20a], and the paper [Ber20b], which illustrates applications of a combined auction/rollout algorithm for solution of multidimensional assignment problems. For example, in an on-line rollout algorithm, at each encountered state, we minimize over the first stage of a multistage dynamic programming problem, and treat the remaining stages approximately, through a relatively fast heuristic. The AWPC algorithm could be a suitable heuristic within this context. 

In conclusion, there are several potentially fruitful possibilities to mesh the AWPC algorithm within the RL methodology. The key property  is that the AWPC algorithm will produce a feasible path starting with arbitrary prices. This path will be near optimal if the starting prices are close to the true (unknown) shortest distances or if they satisfy an $\e$-CS condition with $\e$ relatively small. Moreover, it is plausible that better paths can be obtained by more closely approximating the shortest distances using heuristics and training with data. This conjecture is supported by experience with related auction algorithms, but remains to be established empirically.

\old{
\section{Comparison with Earlier Auction Algorithms for Shortest  Paths}
\pn As noted earlier the AWPC algorithm does not necessarily produce a shortest path. On the other hand, there is an auction algorithm that computes a shortest path and  bears considerable similarity to the AWPC algorithm given here (it is essentially the AWPC algorithm with $\e=0$). We refer to this algorithm as the {\it auction/shortest path} algorithm (ASP for short). The algorithm was first given in the author's paper [Ber91], and is discussed in some detail in the book [Ber98], Section 2.6. In this section, we will compare it to the AWPC algorithm of the present paper, and we will argue that the AWPC algorithm has several advantages, notably faster termination and flexibility to use arbitrary initial prices, at the expense of producing a suboptimal path. 
\old{Then in Section 4.2, we will give an alternative shortest path approach, which is based on converting the shortest path problem to an assignment (weighted matching) problem, which can be solved by using the auction algorithm.
\subsection{Auction Algorithms for a Shortest Path in a Weighted Graph}
} 
The ASP algorithm bears similarity to the AWPC algorithm of Section 3, but does not use an $\e$ parameter, and thus is not suitable for the use of $\e$-scaling. It maintains a single path, starting at the origin, which at each iteration is either extended by adding a new node, or contracted by deleting its terminal node. When the destination becomes  the terminal node of the path, the algorithm terminates. The extensions and contractions, however, are designed to maintain a complementary slackness (CS) property rather than a downhill path property. This property is just the $\e$-CS condition of Section 3.3, with $\e=0$. The CS property is a classical optimality condition in linear programming, and guarantees optimality of the final path, rather than just algorithm termination. Still, however, the AWPC and ASP algorithms produce similar paths, particularly when the parameter $\e$ is relatively small, as we will elaborate later.
The ASP algorithm
maintains a path $P=(s,n_1,\ldots,n_k)$ with no
cycles, and at each iteration modifies $P$ using an extension or a 
contraction operation. The algorithm also maintains a price $p_i$ for each node $i$.
The prices, together with $P$, satisfy at all times the
following CS property 
$$p_i\le a_{ij}+p_j,\qquad \hbox{ for all arcs }(i,j),\xdef\csao{\lab}\eqnum\show{csao}$$  
$$p_i=a_{ij}+p_j,\qquad
\hbox{\rm for all arcs $(i,j)$ of
$P$.}\xdef\csaob{\lab}\eqnum\show{csaob}$$
Note that when CS holds, {\it all the arcs of $P$ are level\/}. Moreover, all the arc discrepancies given by Eq.\ \discrepancy\ are equal to 0. Thus {\it no positive arc discrepancies arise in the ASP algorithm\/}. This is a major difference from the AWPC algorithm, which allows arc discrepancies and the attendant loss of optimality, but also allows more flexible initial prices.
For the ASP algorithm to be valid, we also need to assume that {\it all cycles have positive length\/}, i.e.\ that for every cycle $(i,n_1,\ldots,n_k,i)$ we have 
$$a_{i,n_1}+a_{n_1n_2}+\cdots+a_{n_{k-1}n_k}+a_{n_ki}>0.\xdef\poscycles{\lab}\eqnum\show{csao}$$
This is slightly more restrictive than the nonnegative cycle condition imposed on the AWPC algorithm. The assumption is needed to guarantee that the algorithm does not create a cycle following an extension, since that arcs of $P$ are always level. 
It can be shown that if a pair
$(P,p)$ satisfies the CS conditions, then the portion of $P$ between node $s$ and
any node $i\in P$ is a shortest path from $s$ to $i$, while $p_s-p_i$ is the
corresponding shortest distance. To see this, note that by Eq.\ \csaob,
$p_i-p_k$ is the length of the portion of $P$ between $i$ and $k$,  and that every
path connecting $i$ to $k$ must have length that is greater or equal to $p_i-p_k$ [add 
Eq.\ \csao\ along the arcs of the path].
The ASP algorithm proceeds in iterations, transforming a
pair $(P,p)$ satisfying CS into another pair satisfying CS. At each iteration, the path $P$
is either extended  by a new node or is contracted by deleting its terminal node. In
the latter case the price of the terminal node is increased strictly. A degenerate case
occurs when the path consists of just the origin node $s$; in this case the path
is either extended or is left unchanged with the price $p_s$ being strictly
increased. The iteration (as presented in [Ber98], Section 2.6)  is as follows.
\texshopbox{\pn{\bf Iteration of the Auction/Shortest Path 
Algorithm}
\medskip
\pn Let
$i$ be the terminal node of $P$. If  
$$p_i<\min_{\{j\mid (i,j)\in{\cal
A}\}}\bl\{a_{ij}+p_j\br\},\old{\eqnum\show{csat}}$$  
\pn go to Step 1; else go to Step 2.%}\texshopboxnt{\pn
\smskip
\pn{\bf Step 1 (Contract path):} Set
$$p_i:=\min_{\{j\mid (i,j)\in{\cal A}\}}
\bl\{a_{ij}+p_j\br\},\old{\eqnum\show{csath}}$$ 
and if $i\ne s$, contract $P$. 
Go to the next  iteration. 
\smskip%}\texshopboxnt{\pn
\pn{\bf Step 2 (Extend path):} Extend $P$ by node $j_i$ where
$$j_i=\arg\min_{\{j\mid (i,j) \in{\cal A}\}}
\bl\{a_{ij}+p_j\br\}\old{\eqnum\show{csaf}}$$
(ties are broken arbitrarily). If $j_i$ is the destination $t$, stop; $P$ is the desired
shortest path. Otherwise, go to the next iteration. 
}
\xdef\figureasp{\figr}\figrnum\show{myfigure}
Figure \figureasp\ illustrates the
algorithm. It can be seen from the example of this figure that the terminal node
traces the tree of shortest paths from the origin to the nodes that
are closer to the origin than the given destination.
\topinsert
\centerline{\hskip0pc\epsfxsize = 5.3in \epsfbox{cycle3.eps}}
%\centerline{\tshade{\BoxedEPSF{NN.2.8.eps scaled 750}}}
\vskip0pt
\eightpoint
\def\tablerule{\noalign{\hrule}}
$$\vbox{\offinterlineskip
\hrule
\halign{\vrule\hfill \ #\ \hfill &\vrule\hfill \ #\ \hfill 
&\vrule\hfill \ #\ \hfill &\vrule\hfill 
\ #\ \hfill \vrule\cr
&&\cr
{\bf Iteration \#\lower2ex\hbox{\ }\raise4ex\hbox{\ }}&\hbox{\bf
Path $P$ prior}&\hbox{\bf Price vector $(p_s,p_1,p_2,p_t)$} &\hbox{\bf Type of iteration}\cr
{\bf \lower2ex\hbox{\ }}&\hbox{\bf
to iteration}&\hbox{\bf  prior to iteration} &\cr
\tablerule\cr &&\cr
{1\lower2ex\hbox{\ 
}\raise4ex\hbox{\ }\hbox{\ }}&$(s)$&$(\underline{0},0,0,0)$&Extension to 1\cr 
{2\lower2ex\hbox{\
}\raise2ex\hbox{\ }\hbox{\ }}&$(s,1)$&$(1,\underline{0},0,0)$&Contraction to $s$\cr 
{3\lower2ex\hbox{\
}\raise2ex\hbox{\ }\hbox{\ }}&$(s)$&$(\underline{1},3,0,0)$&Extension to 2\cr 
{4\lower2ex\hbox{\
}\raise2ex\hbox{\ }\hbox{\ }}&$(s,2)$&$(2,3,\underline{0},0)$&Contraction to $s$\cr 
{5\lower2ex\hbox{\
}\raise2ex\hbox{\ }\hbox{\ }}&$(s)$&$(\underline{2},3,2.5,0)$&Extension to 1\cr
{6\lower2ex\hbox{\
}\raise2ex\hbox{\ }\hbox{\ }}&$(s,1)$&$(4,\underline{3},2.5,0)$&Extension to $t$\cr
{7\lower2ex\hbox{\
}\raise2ex\hbox{\ }\hbox{\ }}&$(s,1,t)$&$(4,3,2.5,\underline{0})$&Termination\cr
&&\cr \tablerule\cr}}$$
\fig{-1pc}{\figureasp:}{The results of applying the ASP algorithm to the problem of Fig.\ \figureawsp. The paths and prices generated correspond to $\e=0$ in Fig.\ \figureawsp, which gives the results of applying the AWPC algorithm to the same problem.}\endinsert
It is easily seen that the algorithm maintains CS. Furthermore, the addition of the node  $j_i$
to $P$ following an extension does not create a cycle, since otherwise, in view of the condition 
$p_i\le a_{ij}+p_j$, for every arc $(i,j)$ of the cycle we would have $p_i=a_{ij}+p_j$. By adding
this equality along the cycle, we see that the length of the cycle must be zero, which is not
possible by our assumptions. 
The following proposition, proved in [Ber91] and in [Ber98], Section 2.6, establishes the optimality of the preceding auction algorithm.
\xdef\propauctionproof{\propn}\propnum\show{myproposition}
\texshopbox{\proposition{\propauctionproof:} If there exists at least one path from
the origin to the destination, the preceding auction algorithm terminates
with a shortest path from the origin to the destination.
Otherwise the algorithm never terminates and we have $p_i\to\infty$ for all nodes $i$ in a subset ${\cal N}_\infty$ that contains $s$.}
Thus the ASP algorithm of this section guarantees that the path obtained is shortest with respect to the given arc lengths. However,  the requirement that the initial prices satisfy the CS conditions \csao\ and \csaob\ is a significant restriction in contexts where the arc lengths change over time, or when we attempt to learn good  initial prices with a training process that uses data. 
An interesting observation, which follows from a formal comparison of the respective contraction and extension operations is that {\it the ASP algorithm can be viewed as a limiting form of the AWPC algorithm as $\e\to0$\/}. Indeed let us assume that all cycles have positive length [cf.\ Eq.\ \poscycles], and consider a slightly modified version of the AWPC algorithm whereby when the equality 
$$p_{\hbox{pred}(n_k)}= a_{\hbox{pred}(n_k)n_k}+a_{n_k\hbox{succ}(n_k)}+p_{\hbox{succ}(n_k)}$$
holds [cf.\  case (c2) of the algorithm's statement] we perform an extension rather than a contraction. Suppose that we operate this algorithm formally with $\e=0$, starting with prices that have the CS property \csao-\csaob. It can then be verified that we will obtain the exact same iterates as the ASP algorithm. More generally, it can be verified that if we start the AWPC algorithm with prices that have the CS property  \csao-\csaob, and execute its iterations with $\e$ sufficiently small, we will obtain the same final path as the ASP algorithm, and final prices that exceed those obtained from the ASP algorithm by amounts that are proportional to $\e$. This can also be observed by comparing Figs.\ \figureawsp\ and \figureasp. It is important to emphasize, however, that the AWPC algorithm allows arbitrary initial prices, while the initial prices in the ASP algorithm are restricted by the CS property.
From the preceding discussion, it follows that the optimality of the solution of the ASP algorithm is attained at the expense of slower convergence (as well as  of the restrictive CS requirement for the initial prices). The slow convergence of the ASP algorithm can be overcome by modifications, which are discussed in the paper [BPS95], but the CS requirement for the initial prices remains an obstacle in operating the ASP algorithm in a time-varying environment where on-line replanning is required. 
}
\old{\subsection{Shortest Path Construction by Conversion to an Assignment Problem}
\pn In this section, we show how a shortest path problem can be converted to an assignment problem, which can be solved by the corresponding auction algorithm.
}

\vskip-1pc
\section{Variants of the Algorithms}

\pn In this section we outline variants of the APC and AWPC algorithms of this paper. Detailed development of some of these variants as well as modifications aimed at enhancing computational efficiency will be provided in future reports.

\subsection{Multiple Destinations or Multiple Origins}

\pn The generalization of our algorithms to handle a single origin but multiple destinations is straightforward. We simply maintain a list of destinations that have not yet been reached by the path $P$, and we run the algorithm as if there was a single destination.  Once another destination is reached by $P$, we remove this destination from the list. We then continue similarly, until all destinations are reached. A similar approach has been used to extend the auction/shortest path algorithm of [Ber91] to the case of multiple destinations.

It is also possible to generalize our algorithms to handle multiple origins but a single destination. We simply run the algorithms one origin at a time, as if there was a single origin, and  continue similarly, until a path has been constructed starting from every origin. With multiple origins, however, there is a time-saving possibility. Suppose that we have constructed a path $P_1$ starting at origin 1 and ending at $t$, and then while  constructing a path $P_2$ that starts at origin 2, we land upon a node $n_1\ne t$ of $P_1$ through an extension. Then we can simply join $P_2$ with the tail portion of $P_1$ that starts at $n_1$ and construct a complete path that starts at node 2 and ends at $t$. This can be repeated for all origins, thus eventually constructing a tree of paths to the destination.

\subsection{Forward/Reverse Path Construction Algorithms}

\pn This variant is inspired by forward/reverse versions of the auction algorithm for the assignment problem, due to 
Bertsekas, Casta\~ non, and Tsaknakis [BCT93], which have  also been described and extended in the book  [Ber98]. When applied to the shortest path context, this idea  involves maintaining a forward path that starts from the origin, as well as a reverse path that ends at the destination. The forward path construction uses price increases and proceeds from the origin towards the destination, while the reverse path uses price reductions and proceeds backwards from the destination towards the origin. The forward and the reverse algorithms are symmetric replicas of each other. For the case of a single destination, the algorithm terminates when the forward path that starts from the origin meets the backward path that starts from the destination.

 It is well established by computational practice that forward/reverse variants of auction algorithms generally work faster (and often much faster) than the forward or the reverse algorithms operating alone. It is expected that forward/reverse variants of our APC and AWPC algorithms will similarly work much more efficiently than just corresponding forward algorithms. Thus forward/reverse versions of our algorithms are an important research direction to pursue.

\subsection{Distributed Implementations}

\pn Auction algorithms are  well-suited for parallelization, as extensive implementation studies for the assignment and other network flow problems have shown (see e.g., the book by Bertsekas and Tsitsiklis [BeT89], and the papers by Bertsekas and Casta\~ non [BeC91], Bertsekas et al.\ [BCE95], Beraldi, Guerriero, and Musmanno [BGM97], Zavlanos, Spesivtsev, and Pappas [ZSP08], and Naparstek and Leshem [NaL16]). Possibilities for parallel and distributed asynchronous computation similarly arise within our context when there are multiple destinations and/or multiple origins. The idea is to construct multiple paths simultaneously, with shared price use and asynchronous price updating during the path construction process. Such possibilities are an interesting direction for further research, and have been explored for a related type of auction/shortest path algorithm by Polymenakos and Bertsekas [PoB94].

\subsection{Algorithms Based on Transformations to Equivalent Matching or Assignment Problems}

\pn Path construction problems, weighted and unweighted, can be converted to equivalent assignment and unweighted matching problems, respectively, by using well-known transformations (see e.g., the book [Ber98], which also describes several types of other transformations of network optimization problems to assignment problems). Using these transformations, we can apply the auction algorithm for assignment, which has been investigated extensively. The resulting algorithms bear considerable resemblance to the APC and AWPC  algorithms. A comparative evaluation of these alternative algorithms and their specialized variations is an interesting subject for further research.

\subsection{Variants Involving Multiple Node Price Rises} 

\pn We have discussed so far algorithms that involve a price rise at just the terminal node of the path $P$ maintained by the algorithm. However, a simultaneous contraction and attendant price rise at multiple nodes of $P$ may be possible. Of course this must be done in a way to preserve the downhill path property in some modified form, whereby all arcs of $P$ are either level or downhill, with the last arc of $P$ being downhill following an extension. Such simultaneous price rises may be beneficial if they economize in subsequent computation, and can be facilitated by suitable implementation. For example, when extending $P$ from $n_{k}$ to $n_{k+1}$, we may store the  ``second best neighbor" $n'_{k+1}$ and  corresponding ``second best value" 
$$a_{n_kn'_{k+1}}+p_{n'_{k+1}},$$
where $n'_{k+1}$ is the node that minimizes $a_{n_kn}+p_{n}$ over all $n\ne n_{k+1}$ such that $(n_k,n)$ is an arc. This information can be used at future iterations to determine efficiently if multiple contractions along $P$ can be performed simultaneously. This and other related implementation ideas have been discussed in the papers [Ber95a] and [Ber95b], and have been incorporated in efficient max-flow and minimum cost flow codes, which are available from the author's web site.

\vskip-1pc

 \section{Path-Based Auction Algorithms for Network Transport}

\xdef\examplematching{\exampl}\examplnum\show{myexample}

\xdef\figexamplematching{\figr}\figrnum\show{myfigure}

\pn In this section we will consider a general single-commodity network optimization problem and  a broad extension of our auction/path construction approach for solving it.  The problem involves a network with  a single source, a single sink, and a given amount of supply to be transported from the source to the sink, while respecting given arc capacities. We will use our path construction algorithms as the basis for a methodology to solve (exactly) the unweighted version of this problem, and (inexactly) the weighted version of the problem. Our methods involve successive path constructions, flow augmentations along the constructed paths, and reuse of prices from one path construction to the next. The methods also involve a positive $\e$ parameter, whose choice embodies a tradeoff between accuracy of solution and speed of convergence, and allows for the use of $\e$-scaling.

\subsection{The Unweighted Version of the Problem}

\pn In the unweighted version of the problem we want to transfer a given amount of
flow from a source node to a sink node in a given network without regard for the cost of the transfer. In particular, we have
a directed graph with set of nodes
${\cal N}$  and set of arcs
${\cal A}$. There are two special nodes, denoted $s$ and $t$, which are called the {\it source} and {\it sink\/}, respectively. We assume
that there are no incoming arcs to the source and no outgoing arcs from the sink. Each arc $(i,j)$ is to carry  a flow
$x_{ij}$ that  must satisfy a constraint of the form
$$0\le x_{ij}\le c_{ij},\qquad \forall\ (i,j)\in {\cal
A}.\old{\xdef\uthi{\lab}\eqnum\show{uthi}}$$
 Here  $c_{ij}$ is either a positive scalar or is equal to $\infty$. It represents the ``capacity" of arc $(i,j)$.
 
We are also given a positive scalar $r$, representing supply to be transported from $s$ to $t$, and we consider the problem
of finding a  flow vector $\{x_{ij}\mid (i,j)\in{\cal A}\}$ that satisfies 
$$\sum_{\{j\mid (i,j)\in {\cal A}\}}x_{ij} - \sum_{\{j\mid
(j,i)\in {\cal A}\}}x_{ji}=0,\qquad \forall\ i\in {\cal N},\
i\ne s,\,t,\xdef\uele{\lab}\eqnum\show{uele}$$ 
$$\sum_{\{j\mid (s,j)\in {\cal A}\}}x_{sj}=\sum_{\{j\mid (j,t)\in {\cal
A}\}}x_{jt}=r,\eqnum\show{utwe}$$
$$0\le x_{ij}\le c_{ij},\qquad \forall\ (i,j)\in {\cal
A}.\xdef\uthi{\lab}\eqnum\show{uthi}$$

This is a feasibility problem, whereby we want to transfer a given amount $r$ of
flow from the  source node $s$ to the sink node $t$, while satisfying the arc capacity constraints.\footnote{\dag}{\ninepoint   Sometimes this problem is called the {\it fixed flow problem\/}, to contrast it to the closely related max-flow problem, where we want to maximize the supply $r$, while maintaining feasibility with respect to arc capacities.}
The problem can be solved with a classical approach, which is based on the idea of successive flow augmentations that start at the source, end at the sink and use paths within the so called reduced graph.\footnote{\ddag}{\ninepoint The reduced graph is obtained from the original graph by deleting all arcs $(i,j)$ for which $x_{ij}=c_{ij}$, and by introducing a new (reversed) arc $(j,i)$ for each arc $(i,j)$ such that $0<x_{ij}$. An example of the reduced graph is provided in the subsequent Example \examplematching\ and Fig.\ \figexamplematching. The reduced graph plays a central role in many network flow contexts, including max-flow, primal-dual, and auction algorithms; see [Ber98], starting with Section 3.3.} To this end we can use the APC algorithm of Section 2.1 to construct the paths for the successive flow augmentations, with reuse of the prices from one flow augmentation to next (i.e., use the final  prices from one flow augmentation as the starting prices of the next flow augmentation). A similar auction algorithm with price reuse has been given for max-flow problems in [Ber95a], and tested extensively on large-scale instances with excellent computational results, and far superior performance over competing codes at the time (a code implementing this algorithm is available from the author's website). 

Let us provide an example of this type of algorithm and an illustration of the reduced graph for the case of a matching problem.

\beginexample{\examplematching\ (Successive Path Constructions for Solving a Matching Problem)}\pn Let us consider an unweighted matching problem,  where we want to assign three persons, denoted  $1$, $2$, $3$ to three objects denoted $1'$, $2'$, $3'$. We can transform the problem to a feasibility problem of the form \uele-\uthi, as shown in the top part of Fig.\  \figexamplematching, with all arcs capacities equal to 1. 

The first iteration of the APC algorithm will find some path from $s$ to $t$ in the original graph, and let's assume for the sake of illustration, that this path is $(s,1,2',t)$, matching person $1$ to object $2'$. The reduced graph is then created by reversing the direction of arcs $(s,1)$, $(2',t)$, and $(1,2')$. The second iteration will find some path from $s$ to $t$ in this reduced graph. Let's assume, for the sake of illustration,  that this path is $(s,2,3',t)$, matching person $2$ to object $3'$. The resulting reduced graph is the one shown in the bottom part of Fig.\ \figexamplematching. The third iteration will find  one of the three possible paths from $s$ to $t$ in this reduced graph, resulting in one of the three matchings indicated in Fig.\ \figexamplematching. The algorithm then terminates. 

Generally, for $n\times n$ matching problems, our algorithm will consist of $n$ path constructions, the first $(n-1)$ of which are followed by a suitable modification of the reduced graph. The final prices from each path construction are used as initial prices for the next path construction.
\endexample

\topinsert
\centerline{\hskip0pc\epsfxsize = 4.0in \epsfbox{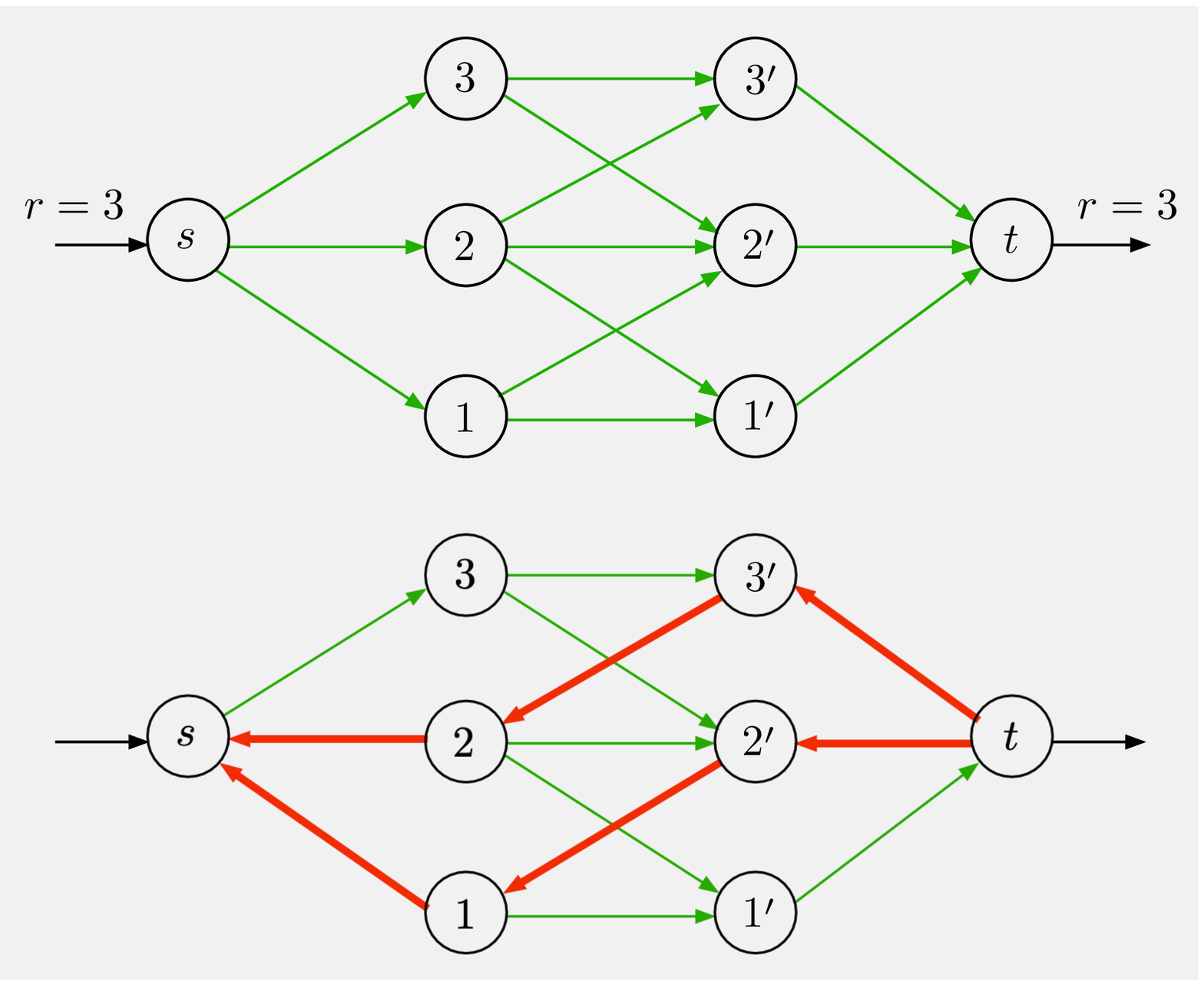}}%\vskip0pt
\vskip0pt
\fig{-0.5pc}{\figexamplematching}{Illustration of the conversion of a $3\times 3$ matching problem to a feasibility problem of the form \uele-\uthi\ (see the top figure and Example \examplematching). All arc capacities are equal to 1 and the supply $r$ is equal to 3. In the bottom part of the figure, we show the reduced graph for the graph shown at the top, after node 1 has been assigned to node 2' and node 2 has been assigned to node 3'. We reverse the direction of the assigned arcs $(1,2')$ and $(2,3')$, and also of the corresponding arcs $(s,1)$, $(s,2)$, $(2',t)$, $(3',t)$. To complete the solution, we need to find a path from $s$ to $t$ in this reduced graph. For any set of initial prices, the APC algorithm will construct one of the following paths:
$$(s,3,3',2, 1',t)\qquad \hbox{resulting in the matching}\qquad (1,2'), (2,1'),(3,3'),$$ 
$$(s,3,3',2, 2',1,1',t)\qquad \hbox{resulting in the matching}\qquad (1,1'), (2,2',(3,3'),$$ 
$$(s,3,2',1, 1',t)\qquad \hbox{resulting in the matching}\qquad (1,1'), (2,3'),(3,2').$$ 
}
\endinsert

\subsection{The Weighted Version of the Problem}

\pn In the linearly weighted version of the problem we have a cost $a_{ij}$ per unit flow on each arc $(i,j)$, and we want to transfer a given amount $r$ of
flow from $s$ to $t$, while minimizing the total cost of the transfer. Thus we want to minimize 
$$\sum_{(i,j)\in{\cal A}}a_{ij}x_{ij},$$
subject to the constraints \uele-\uthi. This is a classical problem of network transport, which contains as special cases problems of assignment (i.e., weighted matching), constrained shortest path (such as $k$ node-disjoint paths), transportation, transhipment, etc, for which specialized auction algorithms have been developed (see [Ber88], [BeC89], [BeC93], and the survey [Ber92]). A well-known primal-dual algorithm to solve the problem is based on successive flow augmentations that start at the source, end at the sink, and use shortest paths within the reduced graph (see [Ber98], Chapter 6). 

To this end we can use the AWPC algorithm of Section 2.1 to construct near-shortest paths for the successive flow augmentations, with reuse of the prices from one flow augmentation to the next. The resulting algorithm will find an approximately optimal solution, with the degree of suboptimality determined by the parameter $\e$ used in the AWPC algorithm. Like AWPC, this algorithm is new, but resembles other already existing auction algorithms. It has its roots in the  auction sequential shortest path algorithms described in the author's paper [Ber95b] (which includes extensive computational comparisons with alternative auction and primal-dual algorithms) and in the book [Ber98], Chapter 7. In particular, the paper [Ber95b] describes several implementation variants and ideas (saving path fragments, early flow augmentations, optimistic extensions, etc). These variants are applicable to the context of the present section and have been incorporated in efficient max-flow and minimum cost flow codes, which are available from the author's web site.

In the case of an assignment problem with a graph such as the one illustrated in Fig.\ \figexamplecycle\ and suitable arc costs $a_{ij}$, the AWPC algorithm will operate similar to the APC algorithm in Example \examplematching. For the $3\times 3$ problem of Fig.\ \figexamplecycle, there will be three augmenting path constructions, guided by the prices and the arc costs. The first two of these will  be followed by a suitable modification of the reduced graph. The final prices from each path construction will be used as initial prices for the next path construction. The price changes involved in the path construction are similar to the ones resulting from the bidding process of the original auction algorithm for the assignment problem [Ber79].

\subsection{Convex Separable Cost Extensions}

\pn We finally note that auction algorithms for the linearly weighted version of the problem have been extended to network problems with convex separable cost functions of the form
$$\sum_{(i,j)\in{\cal A}}f_{ij}(x_{ij}),$$
where $f_{ij}$ is a scalar convex function for each arc $(i,j)\in{\cal A}$. Here we want to find a flow vector $\{x_{ij}\mid (i,j)\in{\cal A}\}$ that minimizes this cost function subject to the constraints \uele-\uthi. We refer to the book [Ber98], Section 9.6 for algorithm descriptions, analysis, and references. These algorithms are based on extended notions of $\e$-complementary slackness and $\e$-scaling, and can be adapted to use path constructions and corresponding flow augmentations, based on the  ideas of this paper.

A further generalized algorithm based on convex separable cost auction ideas addresses flow optimization in networks with gains, and has been given in the paper by Tseng and Bertsekas [TsB00]. This algorithm can also be adapted to use the methodology of the present paper.

\vskip-1pc

\section{Concluding Remarks}
\vskip-0.5pc
\pn In this paper we have introduced an alternative framework for the development of auction algorithms, namely path construction rather than the assignment problem. The new framework allows for arbitrary initial prices, unconstrained by complementary slackness conditions, and is thus well suited for approximations based on training with data, and on-line replanning. The new framework may also be well suited for some application contexts such as shortest path, max-flow, and network transport over long paths. 

Much work remains to be done to explore the possibilities for application of our new auction approach within the broad framework of path planning and network transport. It is also interesting to explore potential applications within reinforcement learning contexts. Moreover, it will be important  to delineate the application areas where the new auction algorithms have superior computational performance over the existing ones, consistent with the empirical observations for max-flow problems [Ber95a]. These objectives are part of the scope of the forthcoming monograph by the author [Ber22a]. 

\vskip-1pc
 
\section{References}
\vskip-0.5pc

\ref[BBS95] Barto, A.\  G., Bradtke, S.\ J., and Singh, S.\ P., 1995.
``Real-Time Learning and Control Using Asynchronous Dynamic Programming," Artificial Intelligence, 
Vol.\ 72, pp.\ 81-138.

\ref[BCE95] Bertsekas, D.\ P., Casta\~ non, D.\ A., Eckstein, J., and Zenios, S.,
1995.\ ``Parallel Computing in Network Optimization," Handbooks in OR and
MS,  Ball, M.\ O., Magnanti, T.\ L., Monma, C.\ L., and Nemhauser, G.\ L.\
(eds.),  Vol.\ 7, North-Holland, Amsterdam, pp.\ 331-399.

\ref[BCT93] 
Bertsekas, D.\ P., Casta\~ non, D.\ A., and Tsaknakis, H.,
1993.\ ``Reverse Auction and the Solution of Inequality Constrained Assignment
Problems," SIAM J.\ on Optimization, Vol.\ 3, pp.\ 268-299.

\ref [BFH03] Brenier, Y., Frisch, U., Henon, M., Loeper, G., Matarrese, S., Mohayaee, R., and Sobolevskii, A., 2003.\ ``Reconstruction of the Early Universe as a Convex Optimization Problem," Monthly Notices of the Royal Astronomical Society, Vol.\ 346, pp.\ 501-524.

\ref[BGM97] Beraldi, P., Guerriero, F., and Musmanno, R., 1997.\ ``Efficient Parallel Algorithms for the Minimum Cost Flow Problem," Journal of Optimization Theory and Applications, Vol.\ 95, pp.\ 501-530.

\ref [BPS95] Bertsekas, D.\ P., Pallottino, S., and Scutell\`a, M.\ G., 1995.\ ``Polynomial
Auction Algorithms for Shortest Paths,'' Computational Optimization and Applications,
Vol.\ 4, pp.\ 99-125.

\ref[BSS08] Bayati, M., Shah, D., and Sharma, M., 2008.\ ``Max-Product for Maximum Weight Matching: Convergence, Correctness, and LP Duality," IEEE Trans.\ on Information Theory, Vol.\ 54, pp.\ 1241-1251.

\ref[BeC89] Bertsekas, D.\ P., and Casta\~ non, D.\ A., 1989.\ ``The Auction Algorithm for the Transportation Problem," Annals of Operations Research, Vol.\ 20, pp.\ 67-96.

\ref[BeC91] Bertsekas, D.\ P., and Casta\~ non, D.\ A., 1991.\ ``Parallel Synchronous and Asynchronous Implementations of the Auction Algorithm," Parallel Computing, Vol.\ 17, pp.\ 707-732.

\ref[BeC93] Bertsekas, D.\ P., and Casta\~ non, D.\ A., 1993.\ ``A Generic Auction Algorithm
for the Minimum Cost Network Flow Problem," Computational Optimization and Applications,
Vol.\ 2, pp.\ 229-260.

\ref[BeE88] Bertsekas, D.\ P., and Eckstein, J., 1988.\ ``Dual Coordinate Step
Methods for Linear Network Flow Problems,'' Math.\ Programming,
Series B, Vol.\ 42, pp.\ 203-243.

\ref [BeT89]  Bertsekas, D.\ P., and Tsitsiklis, J.\ N., 1989.\ Parallel and
Distributed Computation: Numerical Methods, Prentice-Hall, Engl.\ 
Cliffs, N.\ J.\  (can be  downloaded from the author's website).

\ref [BeT97]  
 Bertsimas, D., and Tsitsiklis, J.\ N., 1997.\ Introduction to Linear Optimization,
Athena Scientific, Belmont, MA.

\ref [Ber79] Bertsekas, D.\ P., 1979.\ ``A Distributed Algorithm for the Assignment Problem," Lab. for Information and Decision Systems Report, MIT, May 1979.

\ref [Ber88] Bertsekas, D.\ P., 1988.\ ``The Auction Algorithm: A Distributed Relaxation Method for the Assignment Problem," Annals of Operations Research, Vol.\ 14, pp.\ 105-123.

\ref [Ber91] Bertsekas, D.\ P., 1991.\ ``An Auction Algorithm for Shortest Paths," 
SIAM J.\ on Optimization, Vol.\ 1, pp.\ 425-447.

\ref [Ber92] Bertsekas, D.\ P., 1992.\ ``Auction Algorithms for Network Flow Problems: A Tutorial Introduction," Computational Optimization and Applications, Vol.\ 1, pp.\ 7-66. 

 \ref [Ber95a] Bertsekas, D.\ P., 1995.\ ``An Auction Algorithm for the Max-Flow Problem,"
J.\ of Optimization Theory and Applications, Vol.\ 87, pp.\ 69-101.

 \ref [Ber95b] Bertsekas, D.\ P., 1995.\ ``An Auction/Sequential Shortest Path Algorithm for the Minimum Cost Network Flow Problem," Report LIDS-P-2146, MIT.

\ref [Ber98] Bertsekas, D.\ P., 1998.\
Network Optimization: Continuous and Discrete Models, Athena Scientific,
Belmont, MA (also available on-line from the author's website).

\ref[Ber19] Bertsekas, D.\ P., 2019.\ Reinforcement Learning and Optimal Control, Athena Scientific, Belmont, MA.

\ref[Ber20a] Bertsekas, D.\ P., 2020.\
Rollout, Policy Iteration, and Distributed Reinforcement Learning, Athena Scientific, Belmont, MA.

\ref[Ber20b] Bertsekas, D.\ P., 2020.\ ``Constrained Multiagent Rollout and Multidimensional Assignment with the Auction Algorithm," arXiv preprint, arXiv:2002.07407.

\ref [Ber22a] Bertsekas, D.\ P., 2022.\
Auction Algorithms for Assignment, Path Planning, and Network Transport, Athena Scientific,
Belmont, MA (in preparation).

\ref[Ber22b] Bertsekas, D.\ P., 2022.\
Lessons from AlphaZero for Optimal, Model Predictive, and Adaptive Control, Athena Scientific,
Belmont, MA  (available as an ebook and also on-line from the author's website).

\ref[BiT22] Bicciato, A., and Torsello, A., 2022.\ ``GAMS: Graph Augmentation with Module Swapping," Proc.\ of ICPRAM, pp.\ 249-255.

\ref[CLG22] Clark, A., de Las Casas, D., Guy, A., Mensch, A., Paganini, M., Hoffmann, J., Damoc, B., Hechtman, B., Cai, T., Borgeaud, S,. and Van Den Driessche, G.\ B., 2022. ``Unified Scaling Laws for Routed Language Models," Proc.\ International Conference on Machine Learning, pp.\ 4057-4086.

\ref [Gal16] Galichon, A., 2016.\ Optimal Transport Methods in Economics, Princeton University Press.

\ref[KoY94] Kosowsky, J.\ J., and Yuille, A.\ L., 1994.\ ``The Invisible Hand Algorithm: Solving the Assignment Problem with Statistical Physics," Neural Networks, Vol.\ 7, pp.\ 477-490.

\ref[Kor90] Korf, R.\ E., 1990.\ ``Real-Time Heuristic Search," Artificial Intelligence, Vol.\ 42, pp.\ 189-211.

\ref[LBD21] Lewis, M., Bhosale, S., Dettmers, T., Goyal, N., and Zettlemoyer, L., 2021. ``Base Layers: Simplifying Training of Large, Sparse Models," Proc.\ International Conference on Machine Learning, pp.\ 6265-6274.

\ref[MeT20] Merigot, Q., and Thibert, B., 2021.\ ``Optimal Transport: Discretization and Algorithms," in Handbook of Numerical Analysis,  Elsevier, Vol.\ 22, pp.\ 133-212.

\ref[NaL16] Naparstek, O., and Leshem, A., 2016.\ ``Expected Time Complexity of the Auction Algorithm and the Push Relabel Algorithm for Maximum Bipartite Matching on Random Graphs," Random Structures and Algorithms, Vol.\ 48, pp.\ 384-395.

\ref[PeC19] Peyre, G., and Cuturi, M., 2019.\ Computational Optimal Transport: With Applications to Data Science.\ Foundations and Trends in Machine Learning, Vol.\ 11, pp.\ 355-607.

\ref[PoB94] Polymenakos, L.\ C., and Bertsekas, D.\ P., 1994.\ ``Parallel Shortest Path Auction Algorithms," Parallel Computing, Vol.\ 20, pp.\ 1221-1247.

\ref [Pow11] Powell, W.\ B., 2011.\  Approximate Dynamic Programming: Solving the Curses of Dimensionality, 2nd Edition, J.\ Wiley and Sons, Hoboken, N.\ J.

\ref[San15] Santambrogio, F., 2015.\ Optimal Transport for Applied Mathematicians, Springer Intern.\ Publ.

\ref [Sch16] Schmitzer, B., 2016.\ ``A Sparse Multiscale Algorithm for Dense Optimal Transport," J.\ of Mathematical Imaging and Vision, Vol.\  56, pp.\ 238-259.

\ref [Sch19] Schmitzer, B., 2019.\ ``Stabilized Sparse Scaling Algorithms for Entropy Regularized Transport Problems," SIAM Journal on Scientific Computing, Vol.\ 41, pp.\ A1443-A1481.

\ref [SuB18] Sutton, R., and Barto, A.\  G.,  2018.\ Reinforcement Learning, 2nd Ed., MIT
Press, Cambridge, MA.

\ref[TsB00] Tseng, P., and Bertsekas, D.\ P., 2000.\ ``An $\e$-Relaxation Method for Separable Convex Cost Generalized Network Flow Problems," Mathematical Programming, Vol.\ 88, pp.\ 85-104.

\ref [Vil09] Villani, C., 2009.\ Optimal Transport: Old and New, Springer, Berlin.

\ref [Vil21] Villani, C., 2021.\ Topics in Optimal Transportation, American Mathematical Soc.

\ref[WaD17] Walsh III, J.\ D., and Dieci, L., 2017.\ ``General Auction Method for Real-Valued Optimal Transport," arXiv preprint arXiv:1705.06379.

\ref[WaD19] Walsh III, J.\ D., and Dieci, L., 2019.\ ``A Real-Valued Auction Algorithm for Optimal Transport," Statistical Analysis and Data Mining: The ASA Data Science Journal, 12(6), pp.\ 514-533.

\ref[WaX12] Wang, J., and Xia, Y., 2012.\ ``Fast Graph Construction Using Auction Algorithm," arXiv preprint arXiv:1210.4917.

\ref[ZSP08] Zavlanos, M.\ M., Spesivtsev, L., and Pappas, G.\ J., 2008.\ ``A Distributed Auction Algorithm for the Assignment Problem," Proc.\ 47th IEEE Conference on Decision and Control, pp.\ 1212-1217.

\end

\old{
EXAMPLE OF THE ALGORITHM
\topinsert
\centerline{\hskip0pc\epsfxsize = 4.5in \epsfbox{cycle1.eps}}%\vskip0pt
\vskip-8pt
\eightpoint
\def\tablerule{\noalign{\hrule}}
$$\vbox{\offinterlineskip
\hrule
\halign{\vrule\hfill \ #\ \hfill &\vrule\hfill \ #\ \hfill 
&\vrule\hfill \ #\ \hfill &\vrule\hfill 
\ #\ \hfill \vrule\cr
&&\cr
{\bf Iteration \#\lower2ex\hbox{\ }\raise4ex\hbox{\ }}&\hbox{\bf
Path $P$ prior}&\hbox{\bf Price vector $(p_s,p_1,p_2,p_3,p_t)$} &\hbox{\bf Type of iteration}\cr
{\bf \lower2ex\hbox{\ }}&\hbox{\bf
to iteration}&\hbox{\bf  prior to iteration} &\cr
\tablerule\cr &&\cr
{1\lower2ex\hbox{\ 
}\raise4ex\hbox{\ }\hbox{\ }}&$(s)$&$(\underline{0},0,0,0,0)$&Extension to 1\cr 
{2\lower2ex\hbox{\
}\raise2ex\hbox{\ }\hbox{\ }}&$(s,1)$&$(1+\e,\underline{0},0,0,0)$&Contraction to $s$\cr 
{3\lower2ex\hbox{\
}\raise2ex\hbox{\ }\hbox{\ }}&$(s)$&$(\underline{1+\e},1+\e,0,0,0)$&Extension to 1\cr 
{4\lower2ex\hbox{\
}\raise2ex\hbox{\ }\hbox{\ }}&$(s,1)$&$(2+2\e,\underline{1+\e},0,0,0)$&Extension to 2\cr 
{5\lower2ex\hbox{\
}\raise2ex\hbox{\ }\hbox{\ }}&$(s,1,2)$&$(2+2\e,1+2\e,\underline{0},0,0)$&Contraction to 1\cr
{6\lower2ex\hbox{\
}\raise2ex\hbox{\ }\hbox{\ }}&$(s,1)$&$(2+2\e,\underline{1+2\e},1+\e,0,0)$&Contraction to $s$\cr
{7\lower2ex\hbox{\
}\raise2ex\hbox{\ }\hbox{\ }}&$(s)$&$(\underline{2+2\e},2+2\e,1+\e,0,0)$&Extension to 1\cr
{8\lower2ex\hbox{\
}\raise2ex\hbox{\ }\hbox{\ }}&$(s,1)$&$(3+3\e,\underline{2+2\e},1+\e,0,0)$&Extension to 2\cr
{9\lower2ex\hbox{\
}\raise2ex\hbox{\ }\hbox{\ }}&$(s,1,2)$&$(3+3\e,2+3\e,\underline{1+\e},0,0)$&Extension to 3\cr
{10\lower2ex\hbox{\
}\raise2ex\hbox{\ }\hbox{\ }}&$(s,1,2,3)$&$(3+3\e,2+3\e,1+3\e,\underline{0},0)$&Contraction to 2\cr
{11\lower2ex\hbox{\
}\raise2ex\hbox{\ }\hbox{\ }}&$(s,1,2)$&$(3+3\e,2+3\e,\underline{1+3\e},3+4\e,0)$&Contraction to 1\cr
{12\lower2ex\hbox{\
}\raise2ex\hbox{\ }\hbox{\ }}&$(s,1)$&$(3+3\e,\underline{2+3\e},4+5\e,3+4\e,0)$&Contraction to $s$\cr
{13\lower2ex\hbox{\
}\raise2ex\hbox{\ }\hbox{\ }}&$(s)$&$(\underline{3+3\e},5+6\e,4+5\e,3+4\e,0)$&Extension to 1\cr
{14\lower2ex\hbox{\
}\raise2ex\hbox{\ }\hbox{\ }}&$(s,1)$&$(6+7\e,\underline{5+6\e},4+5\e,3+4\e,0)$&Extension to 2\cr
{$\ldots$\lower2ex\hbox{\
}\raise2ex\hbox{\ }\hbox{\ }}&$\ldots$&$\ldots$&$\ldots$\cr
\tablerule\cr}}$$
\fig{-1pc}{\figexamplecycle}{The shortest path problem of Example \examplecycle\ (top part of the figure). The arc lengths are shown next to the arcs (all lengths are equal to 1, except for the length of arc $(2,t)$ which has a large length $L$. There is only one point where the algorithm can go wrong, at node 2 where there is a choice between going to $t$ or going to 3. The only $s$-to-$t$ path is $(s,1,2,t)$, but if $\e$ is very small, the algorithm explores the possibility of reaching the destination through node 3 for many iterations, while repeating the cycle
$$s\to1\to2\to3\to2\to1\to s\to1\ldots$$
(middle part of the figure). On the other hand, if $3+4\e>L$, then after iteration 14, following an extension to node 2, the algorithm compares the prices of nodes 3 and $t$, performs an extension to $t$, and terminates (top part of the figure).
}
\endinsert
}